\documentclass[11pt,a4paper]{article}
\usepackage{amsmath,amssymb,amsthm,fullpage,graphicx,pstricks}

\newtheorem{thm}{Theorem}[section]
\newtheorem{lem}[thm]{Lemma}
\newtheorem{propn}[thm]{Proposition}
\newtheorem{cor}[thm]{Corollary}

\newtheorem{defn}[thm]{Definition}
\newtheorem{rem}[thm]{Remark}

\newtheorem{assu}[thm]{Assumption}

\let\OLDthebibliography\thebibliography
\renewcommand\thebibliography[1]{
  \OLDthebibliography{#1}
  \setlength{\parskip}{0pt}
  \setlength{\itemsep}{0pt plus 0.3ex}
}

\begin{document}

\title{Scaling limits of stochastic processes\\associated with resistance forms}
\author{D.~A.~Croydon}
\maketitle

\begin{abstract}
\noindent
We establish that if a sequence of spaces equipped with resistance metrics and measures converge with respect to the Gromov-Hausdorff-vague topology, and a certain non-explosion condition is satisfied, then the associated stochastic processes also converge. This result generalises previous work on trees, fractals, and various models of random graphs. We further conjecture that it will be applicable to the random walk on the incipient infinite cluster of critical bond percolation on the high-dimensional integer lattice.
\end{abstract}



\section{Introduction}

In the recent work \cite{ALW} it was shown that if a sequence of tree-like measured metric spaces converge to a limit in the Gromov-Hausdorff-vague topology (see Section \ref{ghvsubsec} for a definition), then so do the associated stochastic processes (cf. \cite[Theorem 6.1]{BCK}). That the convergence of metrics and measures is enough to yield convergence of processes essentially stems from the fact that these provide the natural scale functions and speed measures in this setting. (Hence \cite{ALW} can be viewed as a generalisation of the one-dimensional result of \cite{Stone}.) To put this another way, on the tree-like spaces considered, the distance can be interpreted as a so-called `resistance metric', and this characterises, as the electrical energy, a corresponding quadratic form \cite{AEW} (see also \cite{Kigdendrite}). Together with the specification of an invariant measure, the metric thus determines uniquely a Dirichlet form and stochastic process. From such an observation, one is naturally led to consider whether a similar result to that proved in \cite{ALW} holds for other spaces equipped with a resistance metric and measure. This question was addressed in \cite{CHK}, but only under a strong uniform volume doubling condition (see Remark \ref{explosionrem}(b)). The aim of this article is to remove the latter restriction, and thereby establish the result in much greater generality, enabling a wider range of examples to be handled.

For the setting and notation of the present work, we closely follow \cite{CHK} (for further details, see Section \ref{rfsubsec}). In particular, we let $\mathbb{F}$ be the collection of quadruples of the form $(F,R,\mu,\rho)$, where: $F$ is a non-empty set; $R$ is a resistance metric on $F$ such that closed bounded sets in $(F,R)$ are compact (note this implies $(F,R)$ is complete, separable and locally compact); $\mu$ is a locally finite Borel regular measure of full support on $(F,R)$; and $\rho$ is a marked point in $F$. Note that the resistance metric is associated with a resistance form $(\mathcal{E},\mathcal{F})$ (see Definition \ref{resformdef} below), and we will further assume that for elements of $\mathbb{F}$ this form is regular in the sense of Definition \ref{regulardef}. In particular, this ensures the existence of a related regular Dirichlet form $(\mathcal{E},\mathcal{D})$ on $L^2(F,\mu)$, which we suppose is recurrent, and also a Hunt process $((X_t)_{t\geq 0},\:(P_x)_{x\in F})$.

The intention here is to show that if we have a sequence $(F_n,R_n,\mu_n,\rho_n)_{n\geq 1}$ in $\mathbb{F}$ that converges to some $(F,R,\mu,\rho)\in\mathbb{F}$ with respect to the Gromov-Hausdorff-vague topology, then the associated processes $((X^n_t)_{t\geq 0},\:(P^n_x)_{x\in F_n})$ also converge. The only additional assumption we make, formulated at (\ref{resgrowthcond}) in the following, regards the growth of the resistance from $\rho_n$ to the complement of the open ball $B_n(\rho_n,r)$ in $(F_n,R_n)$, and ensures that the associated processes do not explode (see Remark \ref{explosionrem}(a) for further discussion of this point). It is a natural condition in the context of recurrent processes (see Lemma \ref{reclem}).

{\assu \label{mainassu} The sequence $(F_n,R_n,\mu_n,\rho_n)_{n\geq 1}$ in $\mathbb{F}$ satisfies
\begin{equation}\label{e0}
\left(F_n,R_n,\mu_n,\rho_n\right)\rightarrow\left(F,R,\mu,\rho\right)
\end{equation}
in the Gromov-Hausdorff-vague topology for some $(F,R,\mu,\rho)\in\mathbb{F}$. Moreover, it holds that
\begin{equation}\label{resgrowthcond}
\lim_{r\rightarrow\infty}\limsup_{n\rightarrow\infty}R_n\left(\rho_n,B_n\left(\rho_n,r\right)^c\right)=\infty.
\end{equation}}

We are now ready to state our main result.

{\thm\label{mainthm} Under Assumption \ref{mainassu}, it is possible to isometrically embed $(F_n,R_n)$, ${n\geq 1}$, and $(F,R)$ into a common metric space $(M,d_M)$ in such a way that
\[P^n_{\rho_n}\left(\left(X^n_t\right)_{t\geq 0}\in \cdot\right)\rightarrow P_{\rho}\left(\left(X_t\right)_{t\geq 0}\in \cdot\right)\]
weakly as probability measures on $D(\mathbb{R}_+,M)$ (that is, the space of cadlag processes on $M$, equipped with the usual Skorohod $J_1$-topology).}

\begin{figure}[t]
\begin{center}
\scalebox{0.3}{\includegraphics{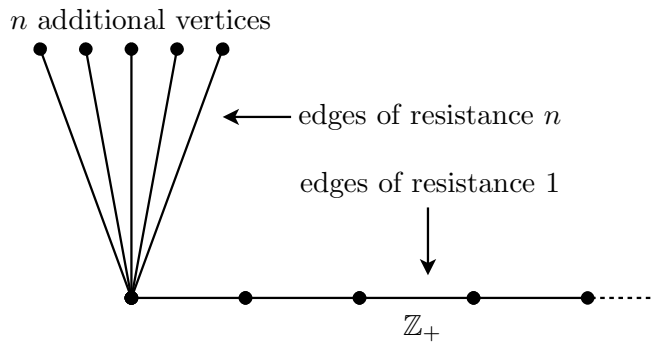}}
\rput(-3.3,-0.3){$\mathbb{Z}_+$}
\rput(-7,3.8){$n$ additional vertices}
\rput(-3.2,2.5){edges of resistance $n$}
\rput(-3.2,1.6){edges of resistance $1$}
\end{center}
\caption{Example of a family of resistance forms failing to satisfy (\ref{resgrowthcond}).}\label{fig1}
\end{figure}

{\rem\label{explosionrem} (a) As we will establish in the proof of Lemma \ref{nrn} below, the resistance growth condition at (\ref{resgrowthcond}) ensures non-explosion of the processes $X^n$ in the limit. To illustrate why it is necessary to have some additional condition beyond merely assuming Gromov-Hausdorff-vague convergence as at (\ref{e0}), consider the sequence of $(F_n,R_n,\mu_n,\rho_n)_{n\geq 1}$ shown in Figure \ref{fig1}. In particular, $F_n$ is $\mathbb{Z}_+$ with $n$ additional vertices, $R_n$ is the resistance on the network obtained by assuming edge resistances as displayed, $\mu_n$ is the measure placing mass one on each vertex, and $\rho_n=0$. Clearly $(F_n,R_n,\mu_n,\rho_n)\rightarrow (F_0,R_0,\mu_0,\rho_0)$ (i.e.\ the limit is $\mathbb{Z}_+$, with no additional vertices). However, the processes $X^n$ do not converge to $X^0$, rather the limiting process has a probability of exploding to infinity of $\frac{1}{2}$ on each visit to the origin. We observe this is consistent with our result, since
\[\lim_{r\rightarrow\infty}\limsup_{n\rightarrow\infty}R_n\left(\rho_n,B_n\left(\rho_n,r\right)^c\right)=1\]
in this case. We also note that in \cite{ALW}, together with Gromov-Hausdorff-vague convergence of trees, a condition (A0) on the lengths of edges leaving compact sets was assumed. Although the condition (A0) of \cite{ALW} did not prevent explosion, it and (\ref{resgrowthcond}) play essentially the same role in the proofs of the relevant results, and neither condition implies the other. As an example of a family of trees satisfying (\ref{resgrowthcond}) but not (A0) in \cite{ALW}, consider $(F_n,R_n,\mu_n,\rho_n)_{n\geq 1}$ as above, but with now only $\sqrt{n}$ additional vertices (still connected by edges of resistance $n$). Evidently the lengths of edges connected to $\rho_n=0$ are unbounded as $n\rightarrow\infty$, and so (A0) in \cite{ALW} does not hold, but (\ref{resgrowthcond}) does hold. Moreover, it is easy to see that the processes $X^n$ do converge to $X^0$ in this case, since the probability of explosion in the limit is zero.\medskip

\noindent
(b) The condition at \eqref{resgrowthcond} is satisfied if all the spaces are compact. Moreover, it is also satisfied if the spaces exhibit uniform volume growth with volume doubling, i.e.\ there exists a function $v(r)$ satisfying $v(2r)\leq c_0v(r)$ and for which $c_1v(r)\leq \mu_n(B_n(x,r))\leq c_2v(r)$, since in that case $c_3r\leq R_n(\rho_n,B_n\left(\rho_n,r\right)^c)\leq r$ (see \cite[Lemma 4.1]{Kum}). The latter condition was the one included in \cite[Assumption 1.2]{CHK}, and so this confirms the assumptions of this paper are weaker.\medskip

\noindent
(c) Whilst the resistance metric is naturally connected to the behaviour of processes, in a number of settings, there might also be another relevant metric. For instance, in studying the random walk on the two-dimensional uniform spanning tree in \cite{BCK}, or on the range of the high-dimensional random walk in \cite{RWRRW}, understanding the resistance metric was crucial, but it was also relevant to consider the embedding into Euclidean space. Similarly, many self-similar fractals are naturally defined as subsets of $\mathbb{R}^d$ (or some other metric space), and it might sometimes be preferable to state results for graphs converging to such fractals in that space, rather than an abstract metric space isometric to the resistance metrics. To allow for this, in Section \ref{embsec} below, we consider an adaptation of the Gromov-Hausdorff-vague topology that includes a continuous (but not necessarily isometric) embedding of our spaces into a common metric space, and prove a convergence result for the embedded processes in this setting (see Theorem \ref{embthm}, cf.\ the result for trees proved as \cite[Theorem 6.1]{BCK}). We also establish a version of this result for random spaces (see Theorem \ref{randembthm}).\medskip

\noindent
(d) As remarked in \cite{Kig}, the theory of resistance forms was introduced in the study of analysis on low-dimensional fractals. Thus it is not surprising that it is within this sphere that our main result finds the greatest applicability. In Section \ref{exsec}, we show how the present work subsumes a number of previous results for trees, fractals, and various models of random graphs. We further present a conjecture about the incipient infinite cluster of bond percolation on the high-dimensional integer lattice.}
\bigskip

To prove Theorem \ref{mainthm}, we follow a similar strategy to \cite{ALW}, that is, we first show tightness of the relevant processes, and then characterise the limit using an occupation density formula. The key to the effectiveness of this approach is that in the setting of resistance forms, as is also the case for trees, we have an explicit expression for the occupation density, or to use alternative terminology, the resolvent kernel of the process killed on hitting a closed set (see Lemma \ref{kernellem} below). Since this is given in terms of the resistance metric, convergence of resistance metrics naturally implies convergence of resolvents (see Proposition \ref{resolvprop}). We note that, in checking this and the tightness of the processes, numerous variations in the technical details are required to the argument of \cite{ALW} to deal with the fact that our spaces are not trees in general. To fill the gaps, we often depend on tools developed for resistance forms in \cite{Kig}. As such our result might also be seen as a generalisation of the results regarding limits of random walks on networks from \cite{Kigdendrite, kig1}, which were motivated by the study of stochastic processes on fractals, and essentially cover the case when $F_n$ is discrete and satisfies $F_n\subseteq F_{n+1}$.

As remarked above, the main result of this paper is moreover related to \cite[Theorem 1.3]{CHK}, which established the same conclusion, but under a uniform volume doubling condition. Whilst the present paper imposes less restrictive assumptions, we note that the advantage of assuming the uniform volume doubling condition is that this also provides equicontinuity and convergence of local times, which our argument does not do. This is potentially relevant in the study of time-changes of processes associated with resistance forms according to time-dependent measures. The argument of \cite{CHK} was also slightly shorter and quite different in character, depending on finite graph approximations, rather than a study of the resolvent.

Finally, we note that the question of whether the distribution of a stochastic process is stable under perturbations in the geometry of the underlying space is also relevant beyond the setting of resistance metrics, and has been considered in various settings. We would particularly like to draw to the attention of the reader the recent work of \cite{Suz2, Suz1}. Indeed, building on \cite{GMS}, it was shown in the aforementioned articles that if a sequence of metric measure spaces satisfying a certain `Riemannian curvature-dimension condition' converge with respect to the measured Gromov-Hausdorff topology, then the associated Brownian motions also converge. Although the statement of the result parallels ours, the techniques applied are again very different, depending on the Mosco convergence of Cheeger energies.

The remainder of the article is organised as follows. In Section \ref{prelimsec}, we present preliminary results on resistance forms and the Gromov-Hausdorff-vague topology. We then, in Section \ref{rksec}, introduce the resolvent kernel, and some of its fundamental properties. In Sections \ref{tightsec} and \ref{semisec}, we establish tightness of the processes $X^n$, $n\geq1$, and show semigroup convergence, in the case when the underlying spaces are compact. Section \ref{proofsec} then ties together the proceeding sections to establish the main result, including extending this to the locally compact setting. The incorporation of embeddings to other spaces is described in Section \ref{embsec}, before finally, in Section \ref{exsec}, we outline some example applications.

\section{Preliminaries}\label{prelimsec}

\subsection{Resistance forms and associated processes}\label{rfsubsec}

Following \cite{CHK}, in this section we recall some basic properties of resistance forms, starting with their definition. The reader is referred to \cite{Kig} for further background.

{\defn [{\cite[Definition 3.1]{Kig}}]\label{resformdef} Let $F$ be a non-empty set. A pair $(\mathcal{E},\mathcal{F})$ is called a \emph{resistance form} on $F$ if it satisfies the following five conditions.
\begin{description}
  \item[RF1] $\mathcal{F}$ is a linear subspace of the collection of functions $\{f:F\rightarrow\mathbb{R}\}$ containing constants, and $\mathcal{E}$ is a non-negative symmetric quadratic form on $\mathcal{F}$ such that $\mathcal{E}(f,f)=0$ if and only if $f$ is constant on $F$.
  \item[RF2] Let $\sim$ be the equivalence relation on $\mathcal{F}$ defined by saying $f\sim g$ if and only if $f-g$ is constant on $F$. Then $(\mathcal{F}/\sim,\mathcal{E})$ is a Hilbert space.
  \item[RF3] If $x\neq y$, then there exists a $f\in \mathcal{F}$ such that $f(x)\neq f(y)$.
  \item[RF4] For any $x,y\in F$,
  \begin{equation}\label{resdef}
  R(x,y):=\sup\left\{\frac{\left|f(x)-f(y)\right|^2}{\mathcal{E}(f,f)}:\:f\in\mathcal{F},\:\mathcal{E}(f,f)>0\right\}<\infty.
  \end{equation}
  \item[RF5] If $\bar{f}:=(f \wedge 1)\vee 0$, then $\bar{f}\in\mathcal{F}$ and $\mathcal{E}(\bar{f},\bar{f})\leq\mathcal{E}({f},{f})$ for any $f\in\mathcal{F}$.
\end{description}}

More generally than (\ref{resdef}), for sets $A,B\subseteq F$, we can write
\[
R(A,B)=\left(\inf\left\{\mathcal{E}(f,f):\:
f\in \mathcal{F},\: f|_A=1,\: f|_B=0\right\}\right)^{-1},
\]
which is the effective resistance between $A$ and
$B$ (interpreted to be zero if the infimum is taken over an empty set). With this definition, we clearly have $R(\{x\},\{y\})=R(x,y)$. The latter function gives a metric on $F$ (see \cite[Proposition 3.3]{Kig}); we call this the \emph{resistance metric} associated with $(\mathcal{E},\mathcal{F})$.  As per the definition of $\mathbb{F}$ in the introduction, we will henceforth assume that we have a non-empty set $F$ equipped with a resistance form $(\mathcal{E},\mathcal{F})$, and denote the corresponding resistance metric $R$. Defining the open ball centred at $x$ and of radius $r$ with respect to the resistance metric by $B_R(x,r):=\left\{y\in F:\:R(x,y)<r\right\}$, and writing its closure as $\bar{B}_R(x,r)$, we will also assume that $\bar{B}_R(x,r)$ is compact for any $x\in F$ and $r>0$ (again noting that this implies $(F,R)$ is complete, separable and locally compact). Furthermore, we will restrict our attention to resistance forms that are regular, as described by the following.

{\defn [{\cite[Definition 6.2]{Kig}}]\label{regulardef} Let $C_0(F)$ be the collection of compactly supported, continuous (with respect to $R$) functions on $F$, and $\|\cdot\|_F$ be the supremum norm for functions on $F$. A resistance form $(\mathcal{E},\mathcal{F})$ on $F$ is called \emph{regular} if and only if $\mathcal{F}\cap C_0(F)$ is dense in $C_0(F)$ with respect to $\|\cdot\|_F$.}
\bigskip

We next introduce related Dirichlet forms and stochastic processes. First, suppose $\mu$ is a Borel regular measure on $(F,R)$ such that $0<\mu(B_R(x,r))<\infty$ for all $x\in F$ and $r>0$. Moreover, write $\mathcal{D}$ to be the closure of $\mathcal{F}\cap C_0(F)$ with respect to the inner product $\mathcal{E}_1$ on $\mathcal{F}\cap L^2(F,\mu)$ given by
$\mathcal{E}_1(f,g):=\mathcal{E}(f,g)+\int_Ffgd\mu$. Under the assumption that $(\mathcal{E},\mathcal{F})$ is regular, we then have from {\cite[Theorem 9.4]{Kig}} that $(\mathcal{E},\mathcal{D})$ is a regular Dirichlet form on $L^2(F,\mu)$ (see \cite{FOT} for the definition of a regular Dirichlet form). Moreover, standard theory gives us the existence of an associated Hunt process $((X_t)_{t\geq 0},\:(P_x)_{x\in F})$ (e.g. \cite[Theorem 7.2.1]{FOT}). Note that such a process is, in general, only specified uniquely for starting points outside a set of zero capacity. However, in this setting every point has strictly positive capacity (see \cite[Theorem 9.9]{Kig}), and so the process is defined uniquely everywhere.

For $(F,R,\mu,\rho)\in \mathbb{F}$, it is an assumption that $(\mathcal{E},\mathcal{D})$ is recurrent. We recall that a Dirichlet form being recurrent is equivalent to $1\in \mathcal{D}_e$ and $\mathcal{E}(1,1)=0$ both holding, where $\mathcal{D}_e$ is the extended Dirichlet space, given by the collection of $\mu$-measurable functions $f:F\rightarrow \mathbb{R}$ such that $|f|<\infty$, $\mu$-a.e., and there exists an $\mathcal{E}$-Cauchy sequence $(f_n)_{n\geq0}$ in $\mathcal{D}$ with $f_n(x) \rightarrow f(x)$, $\mu$-a.e.\ (see \cite[Theorem 1.6.3]{FOT}). The following result yields that, in the present setting, recurrence is also equivalent to the `resistance to infinity' being infinite, i.e.\ \eqref{resgrowthlim} holding. In particular, the latter condition holds whenever $(F,R,\mu,\rho)\in \mathbb{F}$, and the resistance growth condition at (\ref{resgrowthcond}) requires that it holds uniformly for all $(F_n,R_n,\mu_n,\rho_n)_{n\geq 1}$ in the relevant sequence.

{\lem\label{reclem} $(\mathcal{E},\mathcal{D})$ is recurrent if and only if
\begin{equation}\label{resgrowthlim}
\lim_{r\rightarrow\infty}R\left(\rho,B_R\left(\rho,r\right)^c\right)=\infty.
\end{equation}}
\begin{proof} Suppose $(\mathcal{E},\mathcal{D})$ is recurrent, so $1\in \mathcal{D}_e$ and $\mathcal{E}(1,1)=0$. By the definition of $\mathcal{D}_e$ and the regularity of $(\mathcal{E},\mathcal{D})$, it follows that there exists a sequence $(f_n)_{n\geq 1}$ in $C_0(F)\cap\mathcal{D}$ with $\mathcal{E}(f_n,f_n)\rightarrow 0$ and $f_n(x) \rightarrow 1$, $\mu$-a.e. It is further known that, when $(\mathcal{E},\mathcal{D})$ is recurrent, $\mathcal{D}_e=\mathcal{F}$ (see \cite[Proposition 2.13]{KL}), and so we may assume that $(f_n)_{n\geq 1}$ in $C_0(F)\cap\mathcal{F}$. Now, let $x$ be a point with $f_n(x)\rightarrow 1$. From the definition of $R$, we have that $|f_n(\rho)-f_n(x)|^2\leq R(\rho,x)\mathcal{E}(f_n,f_n)\rightarrow 0$, and so it is also the case that $f_n(\rho)\rightarrow 1$. Define $\tilde{f}_n(\cdot)=f_n(\cdot)/f_n(\rho)$, which is in $C_0(F)\cap\mathcal{F}$ (at least for large $n$). Moreover, if the support of $f_n$ is contained in $B_R(\rho,r)$, then it holds that $R(\rho,B_R(\rho,r)^c)^{-1}\leq \mathcal{E}(\tilde{f}_n,\tilde{f}_n)$, which implies
\[\limsup_{r\rightarrow\infty}R\left(\rho,B_R\left(\rho,r\right)^c\right)^{-1}\leq \frac{\mathcal{E}({f}_n,{f}_n)}{f_n(\rho)^2}.\]
Since the upper bound here converges to zero, we obtain (\ref{resgrowthlim}).

To prove the converse, we start by observing (\ref{resgrowthlim}) implies that, for each $r$, there exists a function $f_r\in \mathcal{F}$ such that $f_r(\rho)=1$ and $f|_{B_R(\rho,r)^c}=0$, and also $\mathcal{E}(f_r-1,f_r-1)=\mathcal{E}(f_r,f_r)\rightarrow 0$ as $r\rightarrow\infty$. For such functions, we have that $f_r\in C_0(F)\cap\mathcal{F}\subseteq\mathcal{D}$, and
\[\left|1-f_r(x)\right|^2=\left|f_r(\rho)-f_r(x)\right|^2\leq R(x,\rho)\mathcal{E}(f_r,f_r)\rightarrow 0,\qquad \forall x\in F.\]
Hence $1\in \mathcal{D}_e$ (and $\mathcal{E}(1,1)=0$), as required for $(\mathcal{E},\mathcal{D})$ to be recurrent.
\end{proof}

In the next two lemmas, we recall two further properties of the process $X$. Lemma \ref{commutetimeidentity} is a version of the commute time identity, which is well known for random walks on graphs, and Lemma \ref{ltlem} gives the existence of corresponding local times. Throughout the article, we write
\begin{equation}\label{hitting}
\sigma_A:=\inf\left\{t>0:\:X_t\in A\right\}
\end{equation}
for the hitting time of a set $A$ by $X$, and abbreviate $\sigma_x:=\sigma_{\{x\}}$. Moreover, we denote by $\mathbb{F}_c$ the subset of $\mathbb{F}$ containing elements $(F,R,\mu,\rho)$ such that $(F,R)$ is compact.

{\lem[{\cite[Proof of Lemma 2.9]{CHK}}]\label{commutetimeidentity} If $(F,R,\mu,\rho)\in \mathbb{F}_c$, then
\[E_x\sigma_y+E_y\sigma_x=R(x,y)\mu(F),\qquad \forall x,y\in F.\]}

{\lem[{\cite[Lemma 2.4(b)]{CHK}}]\label{ltlem} If $(F,R,\mu,\rho)\in \mathbb{F}$, then $X$ admits jointly measurable local times $(L_t(x))_{x\in F, t\geq0}$ that satisfy $P_x$-a.s.\ for any $x$,
\[\int_0^t\mathbf{1}_A(X_s)ds=\int_AL_t(y)\mu(dy)\]
for all measurable subsets $A\subseteq F$ and $t\geq 0$.}
\bigskip

Finally, we characterise the process $X$ observed on a compact subset $A\subseteq F$, with respect to a (finite) Borel measure $\nu$ that has $A$ as its support. First, define a continuous additive functional $(\mathcal{A}_t)_{t\geq 0}$ by setting $\mathcal{A}_t:=\int_AL_t(x)\nu(dx)$, and let $(\tau(t))_{t\geq 0}$ be its right-continuous inverse, i.e.\ $\tau(t):=\inf\left\{s>0:\:\mathcal{A}_s>t\right\}$. It is then the case that $(({X}^{\nu}_t)_{t\geq 0},\:(P_x)_{x\in A})$ is also a strong Markov process, where ${X}^{\nu}_t:=X_{\tau(t)}$; this is the \emph{trace} of $X$ on $A$, with respect to $\nu$. Conveniently, we have the following connection, which is an application of the trace theorem for Dirichlet forms (see \cite[Section 6.2]{FOT} for further background on the traces of Dirichlet forms).

{\lem[{\cite[Theorem 2.5 and Lemma 2.6]{CHK}}]\label{restriction} If $(F,R,\mu,\rho)\in \mathbb{F}$, $A\subseteq F$ is a compact set containing $\rho$, and $\nu$ is a finite Borel measure with support $A$, then $(A,R|_{A\times A},\nu,\rho)\in \mathbb{F}_c$, and the associated process is $X^{\nu}$.}

\subsection{Gromov-Hausdorff-vague topology}\label{ghvsubsec}

In this section we introduce the Gromov-Hausdorff-vague topology. For further details, see \cite{ALWGap}. We start by defining a topology on $\mathbb{F}_c$, which we recall is the subset of $\mathbb{F}$ containing elements $(F,R,\mu,\rho)$ such that $(F,R)$ is compact. In particular, for $(F,R,\mu,\rho),(F',R',\mu',\rho')\in\mathbb{F}_c$, we set $\Delta_c((F,R,\mu,\rho),(F',R',\mu',\rho'))$ to be equal to
\begin{equation}\label{deltacdef}
\inf_{M,\psi,\psi'}\left\{d_M^H\left(\psi(F),\psi'(F)\right)+d_M^P\left(\mu\circ\psi^{-1},\mu'\circ\psi'^{-1}\right)+d_M(\psi(\rho),\psi'(\rho'))\right\},
\end{equation}
where the infimum is taken over all metric spaces $M=(M,d_M)$ and isometric embeddings $\psi:(F,R)\rightarrow (M,d_M)$, $\psi':(F',R')\rightarrow (M,d_M)$, and we define $d_M^H$ to be the Hausdorff distance between compact subsets of $M$, and $d_M^P$ to be the Prohorov distance between finite Borel measures on $M$. It is known that $\Delta_c$ defines a metric on the equivalence classes of $\mathbb{F}_c$ (where we say two elements of $\mathbb{F}_c$ are equivalent if there is a measure and root preserving isometry between them).

To extend this topology to one on the equivalence classes of $\mathbb{F}$, we consider bounded restrictions of elements of $\mathbb{F}$. More precisely, for $(F,R,\mu,\rho)\in \mathbb{F}$, define $(F^{(r)},R^{(r)},\mu^{(r)},\rho^{(r)})$ by setting: $F^{(r)}$ to be the closure of the ball in $(F,R)$ of radius $r$ centred at $\rho$, i.e.\ $\bar{B}_R(\rho,r)$; $R^{(r)}$ and $\mu^{(r)}$ to be the restriction of $R$ and $\mu$ respectively to  $F^{(r)}$, and $\rho^{(r)}$ to be equal to $\rho$. By assumption, $F^{(r)}$ is a compact subset of $F$, and $\mu^{(r)}$ is a finite measure of full support on $F^{(r)}$. So, we have from Lemma \ref{restriction} that $(F^{(r)},R^{(r)},\mu^{(r)},\rho^{(r)})\in\mathbb{F}_c$. We then say $(F_n,R_n,\mu_n,\rho_n)\rightarrow(F,R,\mu,\rho)$ in the Gromov-Hausdorff-vague topology if and only if
\[\Delta_c\left((F_n^{(r)},R_n^{(r)},\mu^{(r)}_n,\rho^{(r)}_n),(F^{(r)},R^{(r)},\mu^{(r)},\rho^{(r)})\right)\rightarrow 0\]
for Lebesgue almost-every $r\geq 0$, cf.\ \cite[Definition 5.8]{ALWGap}. From \cite[Proposition 5.9]{ALWGap}, we have the following important consequence of convergence in this topology.

{\lem\label{embeddings} Suppose $(F_n,R_n,\mu_n,\rho_n)$, $n\geq 1$, and $(F,R,\mu,\rho)$ are elements of $\mathbb{F}$ such that $(F_n,R_n,\mu_n,\rho_n)\rightarrow(F,R,\mu,\rho)$ in the Gromov-Hausdorff-vague topology. It is then possible to embed $(F_n,R_n)$, $n\geq 1$, and $(F,R)$ isometrically into the same metric space $(M,d_M)$ in such a way that, for Lebesgue almost-every $r\geq 0$,
\begin{equation}\label{embeddingseq}
d_M^H\left(F_n^{(r)},F^{(r)}\right)\rightarrow 0,\qquad d_M^P\left(\mu_n^{(r)},\mu^{(r)}\right)\rightarrow 0,\qquad d_M(\rho_n^{(r)},\rho^{(r)})\rightarrow 0,
\end{equation}
where we have identified the various objects with their embeddings. Moreover, the space $(M,d_M)$ can be chosen so that closed, bounded subsets of it are compact.}

\section{The resolvent operator and kernel}\label{rksec}

In this section, we study the resolvents of killed processes of resistance forms. More specifically, after recalling various basic properties of such from \cite{Kig} in Lemma \ref{kernellem}, we give an approximation result for the resolvent kernel (see Lemma \ref{approxker}), and demonstrate that, if Assumption \ref{mainassu} holds and the spaces are compact, then the associated resolvents converge (see Proposition \ref{resolvprop}).

Let us start by introducing the objects of interest. Given $(F,R,\mu,\rho)\in\mathbb{F}_c$ and non-empty closed subset $A\subseteq F$, define an operator $G_A$ by setting
\begin{equation}\label{resolvdef}
G_Af(y):=E_y\left(\int_0^{\sigma_A}f\left(X_s\right)ds\right)
\end{equation}
for $y\in F$ and measurable $f:F\rightarrow\mathbb{R}_+$, where we recall from (\ref{hitting}) that $\sigma_A$ is the hitting time of $A$; $G_A$ is the resolvent of the process $X$ killed on hitting $A$. For $x\in F$, we will use the abbreviation $G_x:=G_{\{x\}}$. The following lemma gives an explicit description of the corresponding resolvent kernel. Its proof relies on results from \cite{Kig}.

{\lem \label{kernellem} Let $(F,R,\mu,\rho)\in\mathbb{F}_c$ and $A$ be a non-empty closed subset of $F$. For any measurable $f:F\rightarrow\mathbb{R}_+$, it holds that
\begin{equation}\label{resolveq}
G_Af(y)=\int_Fg_A(y,z)f(z)\mu(dz),\qquad \forall y\in F,
\end{equation}
where
\begin{equation}\label{resolvkernel}
g_A(y,z):=\frac{R(y,A)+R(z,A)-R_A(y,z)}{2},
\end{equation}
with
\[R_A(y,z):=\sup\left\{\mathcal{E}(f,f)^{-1}:\:f\in\mathcal{F},\:f(y)=1,\:f(z)=0,\:f|_A\mbox{ constant}\right\}\]
being the resistance between $y$ and $z$ in the network with $A$ `fused'. Note that the function $g_A$ satisfies $0\leq g_A(y,z)=g_A(z,y)\leq g_A(y,y)=R(y,A)$ and
\begin{equation}\label{contofkernel}
\left|g_A(y,z)-g_A(y,w)\right|\leq R(w,z),\qquad \forall y,z,w\in F.
\end{equation}
Furthermore, if $A=\{x\}$, then $g_x:=g_{\{x\}}$ is given by
\begin{equation}\label{resolvkernelpoint}
g_x(y,z):=\frac{R(x,y)+R(x,z)-R(y,z)}{2}.
\end{equation}}

\begin{proof} From \cite[Theorem 4.1]{Kig}, we have that there exists a unique function $g_A:F\times F\rightarrow\mathbb{R}$ such that, for every $y\in F$, $g_A(y,\cdot)\in\mathcal{F}$ and
\begin{equation}\label{repro}
\mathcal{E}\left(g_A(y,\cdot),u\right)=u(y)
\end{equation}
for every $u\in\{f\in \mathcal{F}:\:f|_A=0\}$. By \cite[Theorem 4.3 and Corollary 5.4]{Kig}, this function can be written as \eqref{resolvkernel}. Moreover, by \cite[Theorem 10.10]{Kig}, we have that
\[g_A(y,z)=\int_0^\infty p_t^{F\backslash A}(y,z)dt,\qquad \forall y,z\in F,\]
where $(p_t^{F\backslash A}(y,z))_{y,z\in F,t>0}$ is the transition density of the process $X$ killed on hitting $A$, which exists and is continuous on $F\times F\times (0,\infty)$ by \cite[Theorem 10.4]{Kig}. This readily implies \eqref{resolveq}, as desired. That $g_A$ satisfies $0\leq g_A(y,z)=g_A(z,y)\leq g_A(y,y)=R(y,A)$ and (\ref{contofkernel}) is proved as part of \cite[Theorem 4.1]{Kig}. The result at (\ref{resolvkernelpoint}) follows immediately from (\ref{resolvkernel}) on noting that $R_{\{x\}}(y,z)=R(y,z)$ by definition.
\end{proof}

We now prove our approximation result for the resolvent kernel.

{\lem\label{approxker} If $(F,R,\mu,\rho)\in\mathbb{F}_c$, then
\[g_x(y,z)-2\varepsilon\leq g_{\bar{B}_R(x,\varepsilon)}(y,z)\leq g_x(y,z)\qquad \forall x,y,z\in F,\:\varepsilon>0.\]}
\begin{proof} The upper bound is easy. Indeed, by the definition of the killed resolvent at (\ref{resolvdef}), the fact that $\mu$ has full support, and the continuity of the resolvent kernel from (\ref{contofkernel}), we have that
\begin{eqnarray*}
0&\leq& \frac{G_x\mathbf{1}_{B_R(z,\delta)}(y)-G_{\bar{B}_R(x,\varepsilon)}\mathbf{1}_{B_R(z,\delta)}(y)}{\mu(B_R(z,\delta))}\\
&=&\frac{\int_{B_R(z,\delta)}\left(g_x(y,w)-g_{\bar{B}_R(x,\varepsilon)}(y,w)\right)\mu(dw)}{\mu(B_R(z,\delta))}\\
&\rightarrow&
g_x(y,z)-g_{\bar{B}_R(x,\varepsilon)}(y,z),
\end{eqnarray*}
as $\delta \rightarrow 0$.

For the lower bound, by the fact that $g_x(y,z)\leq \min\{R(x,y),R(x,z)\}$ and the positivity of $g_{\bar{B}_R(x,\varepsilon)}$ (as noted in Lemma \ref{kernellem}), it will suffice to consider the case when $\min\{R(x,y), R(x,z)\}>\varepsilon$. In this case, consider the function $\tilde{g}$ given by
\[\tilde{g}(z):=\max\left\{\frac{g_x(y,z)-\varepsilon}{R(x,y)-\varepsilon},0\right\}.\]
Note that, since $g_x(y,y)=R(x,y)$, $\tilde{g}(y)=1$, and (\ref{contofkernel}) implies that $\tilde{g}(z)=0$ for $z\in \bar{B}_R(x,\varepsilon)$. Moreover, since we have that $g_x(y,\cdot)\in \mathcal{F}$ (see the proof of Lemma \ref{kernellem}), we have from RF1 and RF5 of Definition \ref{resformdef} that  $\tilde{g}\in\mathcal{F}$, and
\[\mathcal{E}\left(\tilde{g},\tilde{g}\right)\leq \frac{\mathcal{E}\left({g}_x(y,\cdot),{g}_x(y,\cdot)\right)}{\left(R(x,y)-\varepsilon\right)^2}=\frac{{g}_x(y,y)}{\left(R(x,y)-\varepsilon\right)^2}= \frac{R(x,y)}{\left(R(x,y)-\varepsilon\right)^2},\]
where the first equality is a consequence of (\ref{repro}), and we again apply that $g_x(y,y)=R(x,y)$. Hence
\begin{equation}
R\left(y,\bar{B}_R(x,\varepsilon)\right)=\sup\left\{\mathcal{E}(f,f)^{-1}:\:f(y)=1,\:f|_{\bar{B}_R(x,\varepsilon)}=0\right\}
\geq \mathcal{E}\left(\tilde{g},\tilde{g}\right)^{-1}
\geq R(x,y)-2\varepsilon.\label{rtoball}
\end{equation}
Similarly, $R(z,\bar{B}_R(x,\varepsilon))\geq R(x,z)-2\varepsilon$. It is further clear that from the definition that $R_{\bar{B}_R(x,\varepsilon)}(y,z)\leq R(y,z)$, and thus the result follows from the expression for the resolvent kernel given at (\ref{resolvkernel}).
\end{proof}

For the final result of the section, we suppose that Assumption \ref{mainassu} holds, and that $(F_n,R_n)$, ${n\geq 1}$, and $(F,R)$ are isometrically embedded into a common metric space $(M,d_M)$ in the way described in Lemma \ref{embeddings}. We denote by $G_{x_n}^n$ the resolvent associated with $(F_n,R_n,\mu_n,\rho_n)$ when the relevant process is killed on hitting $x_n\in F_n$, and $g^n_{x_n}$ the corresponding kernel.

{\propn\label{resolvprop} Suppose Assumption \ref{mainassu} holds, and that $(F_n,R_n)$, ${n\geq 1}$, and $(F,R)$ are compact. If $x_n,y_n\in F_n$ and $x,y\in F$ are such that $d_M(x_n,x)\rightarrow0$, and $d_M(y_n,y)\rightarrow0$, then
\[G_{x_n}^nf(y_n)\rightarrow G_{x}f(y)\]
for any continuous, bounded $f:M\rightarrow \mathbb{R}$.}
\begin{proof} Let $\Lambda:M^3\rightarrow\mathbb{R}$ be defined by setting
$\Lambda(x,y,z)=({d_M(x,y)+d_M(x,z)-d_M(y,z)})/{2}$, and note that $|\Lambda(x,y,z)-\Lambda(x',y',z')|\leq d_M(x,x')+ d_M(y,y')+ d_M(z,z')$. Moreover, by (\ref{resolvkernelpoint}), we have that $\Lambda(x_n,y_n,\cdot)=g^n_{x_n}(y_n,\cdot)$ on $F_n$, and $\Lambda(x,y,\cdot)=g_{x}(y,\cdot)$ on $F$. Hence, for any continuous, bounded $f:M\rightarrow \mathbb{R}$,
\begin{eqnarray*}
{\left|G_{x_n}^nf(y_n)-G_{x}f(y)\right|}&=&\left|\int_M\Lambda(x_n,y_n,z)f(z)\mu_n(dz)-\int_M\Lambda(x,y,z)f(z)\mu(dz)\right|\\
&\leq&\left(d_M(x_n,x)+d_M(y_n,y)\right)\|f\|_\infty\mu_n(F_n)\\
&&+\left|\int_M\Lambda(x,y,z)f(z)\mu_n(dz)-\int_M\Lambda(x,y,z)f(z)\mu(dz)\right|.
\end{eqnarray*}
For the embedding given by Lemma \ref{embeddings}, we may assume that the space $M$ is compact and $\mu_n\rightarrow\mu$ weakly. Hence both terms above converge to zero as $n\rightarrow\infty$, and we are done.
\end{proof}

\section{Tightness of processes}\label{tightsec}

The aim of this section is to prove the tightness of processes under Assumption \ref{mainassu} in the case when the state spaces are compact. That is, we prove the following, for which we again assume the spaces $(F_n,R_n)$, ${n\geq 1}$, and $(F,R)$ are isometrically embedded into a common metric space $(M,d_M)$ in the way described in Lemma \ref{embeddings}.

{\propn \label{tightness} Suppose Assumption \ref{mainassu} holds, and that $(F_n,R_n)$, ${n\geq 1}$, and $(F,R)$ are compact. For any sequence $(x_n)_{n\geq 1}$ with $x_n\in F_n$, the laws of $X^n$ under $P_{x_n}^n$, $n\geq 1$, are tight in $D(\mathbb{R}_+,M)$.}
\bigskip

The next lemma is a key estimate in the proof of Proposition \ref{tightness}. Its proof is modelled on that of \cite[Lemma 4.4]{ALW}, though requires additional technical input to give us appropriate control of local times.

{\lem\label{hittingtime} (a) Suppose $(F,R,\mu,\rho)\in\mathbb{F}_c$. If $x\in F$, $A$ is a non-empty closed subset of $F\backslash\{x\}$, and $\delta\in(0,R(x,A))$, then
\begin{equation}\label{bound1}
P_x\left(\sigma_{A}\leq t\right)\leq 2\left[1-\left(\frac{R(x,A)-\delta}{R(x,A)+\delta}\right)e^{-\frac{2t}{\mu(B_R(x,\delta))(R(x,A)-\delta)}}\right],\qquad \forall t\geq 0.
\end{equation}
(NB. Under the given assumptions $R(x,A)>0$ by \cite[Theorem 4.3 and Corollary 5.4]{Kig}.) In particular, if $x,y\in F$, $\varepsilon\in (0,R(x,y)/2)$, $\delta\in(0,R(x,y)-2\varepsilon)$, then
\begin{equation}\label{bound2}
P_x\left(\sigma_{\bar{B}_R(y,\varepsilon)}\leq t\right)\leq  4\left[
\frac{\delta}{R(x,y)-2\varepsilon}+\frac{t}{\mu(B_R(x,\delta))(R(x,y)-2\varepsilon-\delta)}\right],\qquad \forall t\geq 0.
\end{equation}
(b) Suppose $(F,R,\mu,\rho)\in\mathbb{F}$. If $\delta\in (0,R(x,{B}_R(\rho,r)^c))$, then
\begin{equation}\label{bound3}
P_\rho\left(\sigma_{{B}_R(\rho,r)^c}\leq t\right)\leq  4\left[
\frac{\delta}{R\left(\rho,{B}_R(\rho,r)^c\right)}+
\frac{t}{\mu(B_R(\rho,\delta))\left(R\left(\rho,{B}_R(\rho,r)^c\right)-\delta\right)}\right],\qquad \forall t\geq 0.
\end{equation}
NB. If ${B}_R(\rho,r)^c=\emptyset$, then we interpret the right-hand side as being equal to zero.}
\begin{proof} First note that, by Lemma \ref{ltlem}, we can write
\[\sigma_A=\int_{F}L_{\sigma_A}(z)\mu(dz)\geq \int_{B_R(x,\delta)}L_{\sigma_A}(z)\mu(dz).\]
So, for any $\lambda>0$, $P_x(\sigma_{A}\leq t)\leq P_x(\lambda\mu\left\{z\in B_R(x,\delta):\:L_{\sigma_A}(z)\geq \lambda\right\}\leq t)$. Taking $\lambda = 2t/\mu(B_R(x,\delta))$, this gives
\begin{align}
P_x\left(\sigma_{A}\leq t\right)&\leq
P_x\left(\mu\left\{z\in B_R(x,\delta):\:L_{\sigma_A}(z)\leq \frac{2t}{\mu(B_R(x,\delta))}\right\}\geq\frac{\mu(B_R(x,\delta))}{2}\right)\nonumber\\
&\leq \frac{2}{\mu(B_R(x,\delta))}\int_{B_R(x,\delta)}P_x\left(L_{\sigma_A}(z)\leq \frac{2t}{\mu(B_R(x,\delta))}\right)\mu(dz)\nonumber\\
&= \frac{2}{\mu(B_R(x,\delta))}\int_{B_R(x,\delta)}\left[1-P_{x}\left(\sigma_z\leq \sigma_A\right)P_z\left(L_{\sigma_A}(z)> \frac{2t}{\mu(B_R(x,\delta))}\right)\right]\mu(dz),\label{return}
\end{align}
where for the second inequality we apply Markov's inequality, and the equality is obtained by applying the strong Markov property at $\min\{\sigma_z,\sigma_A\}$. Now, by applying \cite[Theorem 3.6.3 (equation (3.105) in particular)]{MR} to the process $X$ killed on hitting $A$, we have that $E_xL_{\sigma_A}(z)=g_A(x,z)$. Moreover, the strong Markov property for the killed process yields that, under $P_z$, $L_{\sigma_A}(z)$ is an exponential random variable (see argument around \cite[(3.189)]{MR}, for example). From these two observations, we deduce that
\begin{equation}\label{est1}
P_z\left(L_{\sigma_A}(z)> \frac{2t}{\mu(B_R(x,\delta))}\right)=e^{-\frac{2t}{\mu(B_R(x,\delta))R(z,A)}}.
\end{equation}
Since
\[g_A(x,z)=E_xL_{\sigma_A}(z)=P_x\left(\sigma_z\leq \sigma_A\right)E_zL_{\sigma_A}(z)=P_x\left(\sigma_z\leq \sigma_A\right)g_A(z,z),\]
we also have
\begin{equation}\label{est2}
P_{x}\left(\sigma_z\leq \sigma_A\right)=\frac{g_A(x,z)}{g_A(z,z)}=\frac{R(x,A)+R(z,A)-R_A(x,z)}{2R(z,A)}
\end{equation}
(cf. \cite[(3.109)]{MR}). To estimate $R(z,A)$ in (\ref{est1}) and (\ref{est2}), we first note that \cite[Theorem 4.3 and Corollary 5.4]{Kig} implies that $R_A$ is a metric on $(F\backslash A)\cup \{A\}$. Hence, applying the triangle inequality with this metric (and noting that $R_A(x,z)\leq R(x,z)$), we find that, for $z\in B_R(x,\delta)$,
\[R(x,A)=R_A(x,A)\leq R_A(x,z)+R_A(z,A)\leq R(x,z)+R(z,A)<\delta+R(z,A),\]
and similarly $R(z,A)\leq R(x,A)+\delta$ . Hence, we obtain that, for $z\in B_R(x,\delta)$,
\[P_{x}\left(\sigma_z\leq \sigma_A\right)P_z\left(L_{\sigma_A}(z)> \frac{2t}{\mu(B_R(x,\delta))}\right)\geq\left(\frac{R(x,A)-\delta}{R(x,A)+\delta}\right)e^{-\frac{2t}{\mu(B_R(x,\delta))(R(x,A)-\delta)}}.\]
Substituting this bound into (\ref{return}), we obtain (\ref{bound1}).

For (\ref{bound2}), we start with the elementary observation that
\[\frac{R(x,A)-\delta}{R(x,A)+\delta}\geq 1-\frac{2\delta}{R(x,A)},\]
which, together with $1-x\leq e^{-x}$, implies
\begin{equation}\label{e2}
1-\left(\frac{R(x,A)-\delta}{R(x,A)+\delta}\right)e^{-\frac{2t}{\mu(B_R(x,\delta))(R(x,A)-\delta)}}\leq
\frac{2\delta}{R(x,A)}+\frac{2t}{\mu(B_R(x,\delta))(R(x,A)-\delta)}.
\end{equation}
Moreover, from \eqref{rtoball}, we have that $R(x,\bar{B}_R(y,\varepsilon))\geq R(x,y)-2\varepsilon$. Combining \eqref{bound1}, \eqref{e2} and the latter observation yields the result.

For part (b), consider the space $(\tilde{F},\tilde{R},\tilde{\mu},\tilde{\rho})$ obtained by setting: $\tilde{F}:=B_R(\rho,r)\cup\{B_R(\rho,r)^c\}$, $\tilde{R}:=R_{B_R(\rho,r)^c}$, $\tilde{\mu}:=\mu(\cdot\cap B_R(\rho,r))+\delta_{\{B_R(\rho,r)^c\}}$ (where the second term here is the probability measure placing all its mass on the point set $\{B_R(\rho,r)^c\}$), and $\tilde{\rho}:=\rho$. From \cite[Theorems 4.3, 6.3 and 6.4]{Kig}, we have that $(\tilde{F},\tilde{R},\tilde{\mu},\tilde{\rho})\in\mathbb{F}_c$. Hence we can apply (\ref{bound1}) to obtain
\begin{equation}\label{bound4}
\tilde{P}_{\tilde{\rho}}\left(\tilde{\sigma}_{\{{B}_R(\rho,r)^c\}}\leq t\right)\leq  4\left[
\frac{\delta}{R\left(\rho,{B}_R(\rho,r)^c\right)}+
\frac{t}{\mu(B_R(\rho,\delta))\left(R\left(\rho,{B}_R(\rho,r)^c\right)-\delta\right)}\right],
\end{equation}
where we denote by $\tilde{P}_{\tilde{\rho}}$ and $\tilde{\sigma}_{\{{B}_R(\rho,r)^c\}}$ the law and hitting time of the process associated with $(\tilde{F},\tilde{R},\tilde{\mu},\tilde{\rho})$. (For this, note that ${R}(\rho,{B}_R(\rho,r)^c)=\tilde{R}(\tilde{\rho},\{{B}_R(\rho,r)^c\})$, and also that ${\mu}(B_R(\rho,\delta))\leq \tilde{\mu}(B_{\tilde{R}}(\tilde{\rho},\delta))$.) Finally, from \cite[Section 4.4]{FOT}, we have that the domain of the Dirichlet form of $X$ killed on exiting $B_R(\rho,r)$ is given by
\[\mathcal{F}_{kill}:=\{f\in \mathcal{D}:\:f|_{B_R(\rho,r)^c}=0\}=\{f\in \mathcal{F}:\:f|_{B_R(\rho,r)^c}=0\},\]
where the equality holds because $\mathcal{D}=\mathcal{D}_e\cap L^2(F,\mu)=\mathcal{F}\cap L^2(F,\mu)$ (the first equality here is \cite[Theorem 1.5.2(iii)]{FOT}, and the second is a consequence of \cite[Proposition 2.13]{KL}). Similarly, the Dirichlet form of the process associated with $(\tilde{F},\tilde{R},\tilde{\mu},\tilde{\rho})$, is given by $({\mathcal{E}},\tilde{\mathcal{F}})$, where $\tilde{\mathcal{F}}:=\{f\in\mathcal{F}:\:f|_{B_R(\rho,r)^c}\mbox{ constant}\}$, and the corresponding killed process has domain
\[\{f\in \tilde{\mathcal{F}}:\:f|_{B_R(\rho,r)^c}=0\}=\mathcal{F}_{kill}.\]
This implies the Dirichlet forms of the two killed processes are the same, and therefore, up to exiting $B_R(\rho,r)$, the processes associated with $(F,R,\mu,\rho)$ and $(\tilde{F},\tilde{R},\tilde{\mu},\tilde{\rho})$ are also the same. In particular, the left-hand side of (\ref{bound4}) is equal to the left-hand side of (\ref{bound3}), and the result follows.
\end{proof}

We next apply the hitting time bound of Lemma \ref{hittingtime} to control fluctuations in the sample paths of $X$. We write $N(F,\varepsilon)$ to denote the minimal size of an $\varepsilon$-net of $(F,R)$ (where we recall an $\varepsilon$-net of $(F,R)$ is a subset $A\subseteq F$ such that $\cup_{x\in A}B_R(x,\varepsilon)=F$). In the case $(F,R)$ is compact, $N(F,\varepsilon)<\infty$ for every $\varepsilon>0$.

{\lem \label{fluct} Suppose $(F,R,\mu,\rho)\in\mathbb{F}_c$. For any $\varepsilon>0$, $\delta\in (0,\varepsilon/8]$,
\[\sup_{x\in F}P_x\left(\sup_{s\leq t}R(x,X_s)\geq \varepsilon\right)\leq \frac{32 N(F,\varepsilon/4)}{\varepsilon}\left(\delta+\frac{t}{\inf_{x\in F}\mu(B_R(x,\delta))}\right).\]}
\begin{proof} Let $(x_i)_{i=1}^{N(F,\varepsilon/4)}$ be an $(\varepsilon/4)$-net
of $F$. Clearly $B_R(x,\varepsilon)^c\subseteq\cup_{i:\:R(x,x_i)\geq 3\varepsilon/4}\bar{B}_R(x_i,\varepsilon/4)$. Hence
\[P_x\left(\sup_{s\leq t}R(x,X_s)\geq \varepsilon\right)\leq\sum_{i:\:R(x,x_i)\geq 3\varepsilon/4}P_x\left(\sigma_{\bar{B}_R(x_i,\varepsilon/4)}\leq t\right),\]
and applying \eqref{bound2} with $\varepsilon$ replaced by $\varepsilon/4$, we obtain the result.
\end{proof}

The following result extends the previous bound to a uniform one for sequences of compact resistance forms satisfying Assumption \ref{mainassu}.

{\lem \label{fluctuni} Suppose Assumption \ref{mainassu} holds, and that $(F_n,R_n)$, ${n\geq 1}$, and $(F,R)$ are compact. For any $\varepsilon>0$,
\[\lim_{t\rightarrow0}\limsup_{n\rightarrow \infty}\sup_{x\in F_n}P^n_x\left( \sup_{s\leq t}R_n(x,X^n_s)\geq\varepsilon\right)=0.\]}

\begin{proof} From Lemma \ref{fluct}, it will be sufficient to check that, for each $\varepsilon>0$, $\delta\in (0,\varepsilon/8)$,
\[\sup_{n\geq 1}N(F_n,\varepsilon/4)<\infty,\]
and also
\[\inf_{n\geq 1}\inf_{x\in F_n}\mu_n(B_n(x,\delta))>0.\]
However, both these are ready consequences of the convergence in Gromov-Hausdorff-vague topology; see \cite[Proposition 7.4.12]{BBI} and \cite[Corollary 5.7]{ALWGap} respectively.
\end{proof}

We can now complete the proof of the tightness result of Proposition \ref{tightness}.

\begin{proof}[Proof of Proposition \ref{tightness}] It is easy to check Lemma \ref{fluctuni} implies Aldous' tightness criterion, which implies the result (see \cite[Theorem 16.10 and 16.11]{Kall}, for example).
\end{proof}

\section{Convergence of semigroups}\label{semisec}

We now turn to the study of the semigroups associated with a sequence $(F_n,R_n,\mu_n,\rho_n)_{n\geq 1}$ in $\mathbb{F}$ satisfying Assumption \ref{mainassu}. Since each semigroup is defined on a different space, following \cite{ALW}, we capture these via functions on $\mathcal{M}_1(M)$, the space of Borel probability measures on $M$. In particular, assuming all the spaces are embedded into a common metric space in the way described by Lemma \ref{embeddings}, for each $n\geq 1$, we define a map $Q_n:F_n\times\mathbb{R}_+\rightarrow \mathcal{M}_1(M)$ by setting
\[Q_n(x,t):=P_x^n\left(X_t^n\in \cdot\right).\]
Similarly, for the limiting space $(F,R,\mu,\rho)\in\mathbb{F}$, we define $Q:F\times\mathbb{R}_+\rightarrow \mathcal{M}_1(M)$ by setting
\[Q(x,t):=P_x\left(X_t\in \cdot\right).\]
The main aim of this section is to establish the following limit result, where we denote by $d_M^P$ the Prohorov metric on $\mathcal{M}_1(M)$.

{\propn \label{convtoQ} Suppose Assumption \ref{mainassu} holds, and that $(F_n,R_n)$, ${n\geq 1}$, and $(F,R)$ are compact. It is then the case that, for every $T>0$,
\[\lim_{\delta\rightarrow 0}
\limsup_{n\rightarrow\infty}
\sup_{\substack{x\in F,\:y\in F_{n}:\\d_M(x,y)\leq \delta}}
\sup_{\substack{s,t\in[0,T]:\\|s-t|\leq \delta}}
{d_M^P}\left({Q}(x,s),Q_{n}(y,t)\right)=0.\]}

The strategy we adopt for proving the above result will be similar to that of \cite{ALW}, namely we first establish the existence of subsequential limits (Lemma \ref{subseqlem}, cf.\ \cite[Lemmas 5.3 and 5.4]{ALW}), then show any limit points are Feller (Lemma \ref{fellerlem}, cf.\ \cite[Proposition 5.2]{ALW}) and have the same resolvent as $X$ (Lemma \ref{equalresolv}, cf.\ \cite[Proposition 5.6]{ALW}), and finally check these properties uniquely characterise the limit (cf.\ \cite[Proposition 5.1]{ALW}). Towards the first of these goals, the following equicontinuity result is key. Note that in the proof we write $\sigma^n_A$ for the hitting time of a set $A\subseteq F_n$ by $X^n$, and $\sigma_{x}^n:=\sigma_{\{x\}}^n$.

{\lem\label{equilem} Suppose Assumption \ref{mainassu} holds, and that $(F_n,R_n)$, ${n\geq 1}$, and $(F,R)$ are compact. It is then the case that
\[\lim_{\delta\rightarrow 0}\sup_{n\geq 1}\sup_{\substack{x,y\in F_n:\\R_n(x,y)\leq \delta}}\sup_{\substack{s,t\in\mathbb{R}_+:\\|s-t|\leq \delta}}{d_M^P}\left(Q_n(x,s),Q_n(y,t)\right)=0.\]}
\begin{proof} As in the proof of \cite[Lemma 5.3]{ALW}, the proof proceeds via a coupling of $Q_n(x,s)$ and $Q_n(y,t)$. In particular, for $(F,R,\mu,\rho)\in\mathbb{F}_c$, $x,y\in F$, $0\leq s\leq t$, and $\varepsilon,\eta>0$, note that we have
\[P_x\left(R(X_s,X_{t+\sigma_y})\geq \varepsilon \right)\leq
P_x\left(\sigma_y>\eta \right)+P_x\left(\sup_{u\in [t,t+\eta]}R(X_s,X_u)\geq \varepsilon \right).\]
Now, by the commute time identity (Lemma \ref{commutetimeidentity}), the first term here is bounded above by $R(x,y)\mu(F)/\eta$. Moreover, by the Markov property and Lemma \ref{fluct}, the second term is bounded above by
\[\sup_{z\in F}P_z\left(\sup_{u\in [0,t+\eta-s]}R(z,X_u)\geq \varepsilon \right)\leq \frac{32 N(F,\varepsilon/4)}{\varepsilon}\left(\eta'+\frac{t+\eta-s}{\inf_{x\in F}\mu(B_R(x,\eta'))}\right)\]
for any $\eta'\leq \varepsilon/8$. Hence, under Assumption \ref{mainassu}, we obtain (similarly to the proof of Lemma \ref{fluct}) that
\[\sup_{\substack{x,y\in F_n:\\R_n(x,y)\leq\delta}}\sup_{\substack{s,t\geq 0:\\|s-t|\leq \delta}}P^n_x\left(R(X^n_s,X^n_{t+\sigma^n_y})\geq \varepsilon \right)\leq
c_1\eta^{-1}\delta+c_2 \left(\eta'+\frac{\eta+\delta}{\inf_{x\in F_n}\mu_n(B_n(x,\eta'))}\right),\]
for some finite constants $c_1,c_2$ that are independent of $n$, $\delta$, $\eta$ and $\eta'$ (though may depend on $\varepsilon$). Taking $\eta'$ small, and then $\eta=\sqrt{\delta}$ small, the bound above can be made arbitrarily small (uniformly in $n$). Finally, we observe that if $(X^n_t)_{t\geq0}$ has law $P^n_x$, then the process $(X^n_{t+\sigma^n_y})_{t\geq 0}$ has law $P^n_y$ (note that the commute time identity gives that the stopping time is almost-surely finite). Hence the previous bound implies the result.
\end{proof}

Given Lemma \ref{equilem}, it is straightforward to deduce the existence of subsequential limits by a standard Arzel\`{a}-Ascoli type-argument (cf.\ \cite[Lemma 5.4]{ALW}). We therefore present the statement without proof.

{\lem \label{subseqlem} Suppose Assumption \ref{mainassu} holds, and that $(F_n,R_n)$, ${n\geq 1}$, and $(F,R)$ are compact. There exists a subsequence $(n_m)_{m\geq 1}$ and function $\tilde{Q}:F\times \mathbb{R}_+\rightarrow \mathcal{M}_1(M)$ such that, for every $T>0$,
\begin{equation}\label{subseqlimit}
\lim_{\delta\rightarrow 0}
\limsup_{m\rightarrow\infty}
\sup_{\substack{x\in F,\:y\in F_{n_m}:\\d_M(x,y)\leq \delta}}
\sup_{\substack{s,t\in[0,T]:\\|s-t|\leq \delta}}
{d_M^P}\left(\tilde{Q}(x,s),Q_{n_m}(y,t)\right)=0.
\end{equation}}

We will call any function $\tilde{Q}:F\times \mathbb{R}_+\rightarrow \mathcal{M}_1(M)$ that satisfies (\ref{subseqlimit}) a subsequential limit of $(Q_n)_{n\geq 1}$. Note that it is immediate from (\ref{subseqlimit}) that any such function is continuous. Moreover, if $x_{n_m}\in F_{n_m}$, $x\in F$ is such that $d_M(x_{n_m},x)\rightarrow 0$, then $Q_{n_m}(x_{n_m},t)\rightarrow \tilde{Q}(x,t)$ for all $t\geq 0$. For a given $\tilde{Q}:F\times \mathbb{R}_+\rightarrow \mathcal{M}_1(M)$, define a family of operators $(\tilde{S}_t)_{t\geq 0}$ on $C(F)$ (the space of continuous functions on $F$, which we assume is equipped with the supremum norm) by setting
\[\tilde{S}_tf(x):=\int_Ffd\tilde{Q}(x,t)\qquad \forall f\in C(F),\:x\in F,\:t\geq 0.\]
The following lemma shows this definition yields a nice semigroup for any subsequential limit.

{\lem \label{fellerlem} Suppose Assumption \ref{mainassu} holds, and that $(F_n,R_n)$, ${n\geq 1}$, and $(F,R)$ are compact. If $\tilde{Q}$ is a subsequential limit of $(Q_n)_{n\geq 1}$, then $(\tilde{S}_t)_{t\geq 0}$ is a conservative Feller semigroup (see \cite[Chapter 19]{Kall} for background). Moreover, if $(n_m)_{m\geq 1}$ is the subsequence along which convergence takes place, and $x_{n_m}\in F_{n_m}$ is such that $d_M(x_{n_m},x)\rightarrow 0$ for some $x\in F$, then
\begin{equation}\label{felproc}
P^{n_m}_{\rho_{n_m}}\left(\left(X^{n_m}_t\right)_{t\geq 0}\in \cdot\right)\rightarrow \tilde{P}_{\rho}\left(\left(\tilde{X}_t\right)_{t\geq 0}\in \cdot\right)
\end{equation}
weakly as probability measures on $D(\mathbb{R}_+,M)$, where $((\tilde{X}_t)_{t\geq 0},(\tilde{P}_x)_{x\in F})$ is the Markov process associated with the semigroup $(\tilde{S}_t)_{t\geq 0}$.}

\begin{proof} It is clear from the construction that $(\tilde{S}_t)_{t\geq0}$ are a family of positive contraction operators on $C(F)$. Moreover, since each $\tilde{Q}(x,t)$ is a probability measure, we have that $\tilde{S}_t1=1$. Hence to establish that $(\tilde{S}_t)_{t\geq 0}$ is a conservative Feller semigroup, it remains to show the semigroup property:
\begin{equation}\label{semigroup}
\tilde{S}_{t+s}f=\tilde{S}_{t}\tilde{S}_{s}f,\qquad\forall f\in C(F),\:s,t\geq 0,
\end{equation}
and also the two Feller regularity properties:
\begin{equation}\label{f1}
\tilde{S}_tf\in C(F),\qquad\forall f\in C(F),\;t\geq0,
\end{equation}
\begin{equation}\label{f2}
\lim_{t\rightarrow0}\tilde{S}_tf(x)= f(x),\qquad\forall f\in C(F),\: x\in F.
\end{equation}

Recall from the remark above the lemma that $\tilde{Q}$ is continuous. Hence if we have a sequence $(x_n)_{\geq 1}$ in $F$ such that $x_n\rightarrow x\in F$, then $\tilde{Q}(x_n,t)\rightarrow \tilde{Q}(x,t)$ weakly in $\mathcal{M}_1(M)$. It follows immediately that $\tilde{S}_tf(x_n)\rightarrow \tilde{S}_tf(x)$ for any $f\in C(F)$, which confirms (\ref{f1}) holds.

Towards proving (\ref{semigroup}) and (\ref{f2}), we start by checking that, for any $f\in \mathrm{Lip}_1(M)$ (that is, any $1$-Lipschitz function on $M$), we have that
\begin{equation}\label{aaa}
\lim_{\delta\rightarrow 0}
\limsup_{m\rightarrow\infty}
\sup_{\substack{x\in F,\:y\in F_{n_m}:\\d_M(x,y)\leq \delta}}
\sup_{\substack{s,t\in[0,T]:\\|s-t|\leq \delta}}\left|\tilde{S}_sf(x)-S^{n_m}_tf(y)\right|=0,
\end{equation}
where $S^n$ is the semigroup of $X^n$. In particular, note that, for any $f\in \mathrm{Lip}_1(M)$, $x\in F$, $y\in F_{n_m}$, $s,t\geq0$,
\begin{eqnarray*}
\left|\tilde{S}_sf(x)-S^{n_m}_tf(y)\right|&=&\left|\int_Mf d\tilde{Q}(x,s)-\int_Mf d{Q}^{n_m}(y,t)\right|\\
&\leq& \left(1+\mathrm{diam}_{d_M}(M)\right){d_M^P}\left(\tilde{Q}(x,s),{Q}^{n_m}(y,t)\right),
\end{eqnarray*}
where the second inequality is comes from a standard comparison of the Wasserstein and Prohorov metrics (see \cite[Theorem 2]{GSu}, for example), with $\mathrm{diam}_{d_M}(M)$ being the diameter of the metric space $(M,d_M)$. From this bound and Lemma \ref{subseqlem}, we obtain \eqref{aaa}.

We next prove \eqref{f2}. First note that, by the continuity of $\tilde{Q}$, we have that
$\tilde{S}_tf(x)\rightarrow \tilde{S}_0f(x)$ as $t\rightarrow 0$, for any $f\in C(F)$, and $x\in F$. Thus we are required to show $\tilde{S}_0f(x)=f(x)$. Since we can approximate any continuous function on a compact metric space uniformly by Lipschitz functions \cite{Georg}, it will suffice to prove the latter result when $f\in \mathrm{Lip}_1(F)$. Note that, if $f\in  \mathrm{Lip}_1(F)$, then it can be extended to a function in $\mathrm{Lip}_1(M)$. Thus, if $x_{n_m}\in F_{n_m}$, $x\in F$ is such that $d_M(x_{n_m},x)\rightarrow 0$, then from (\ref{aaa}) we obtain \[\tilde{S}_0f(x)=\lim_{m\rightarrow\infty}S^{n_m}_0f(x_{n_m})=\lim_{m\rightarrow\infty}f(x_{n_m})=f(x),\]
as desired.

To complete the proof that $(\tilde{S}_t)_{t\geq 0}$ is Feller, we check the semigroup property (\ref{semigroup}). Again taking $x_{n_m}\in F_{n_m}$, $x\in F$ such that $d_M(x_{n_m},x)\rightarrow 0$, we have from (\ref{aaa}) and the semigroup property of $(S^{n_m}_t)_{t\geq 0}$ that
\[\tilde{S}_{t+s}f(x)=\lim_{m\rightarrow\infty}{S}^{n_m}_{t+s}f(x_{n_m})=
\lim_{m\rightarrow\infty}{S}^{n_m}_{t}{S}^{n_m}_{s}f(x_{n_m})\]
for any $f\in \mathrm{Lip}_1(M)$. Now, since $\tilde{S}_sf\in C(F)$ (by (\ref{f1})), it can be extended to a function in $C(M)$, and we obtain
\begin{eqnarray}
\left|{S}^{n_m}_{t}{S}^{n_m}_{s}f(x_{n_m})-{S}^{n_m}_{t}\tilde{S}_{s}f(x_{n_m})\right|
&\leq&\sup_{y\in F_{n_m}}\left|{S}^{n_m}_{s}f(y)-\tilde{S}_{s}f(y)\right|.\label{tort}
\end{eqnarray}
To bound the right-hand side here, we note that, under Assumption \ref{mainassu}, it is possible to define a function $g_n:F_n\rightarrow F$ such that
\begin{equation}\label{gndef}
\lim_{n\rightarrow\infty}\sup_{y\in F_n}d_M(y,g_n(y))=0.
\end{equation}
Returning to \eqref{tort}, we can write
\begin{eqnarray}
{\left|{S}^{n_m}_{t}{S}^{n_m}_{s}f(x_{n_m})-{S}^{n_m}_{t}\tilde{S}_{s}f(x_{n_m})\right|}
&\leq&\sup_{y\in F_{n_m}}\left|{S}^{n_m}_{s}f(y)-\tilde{S}_{s}f(g_{n_m}(y))\right|\nonumber\\
&&+\sup_{y\in F_{n_m}}\left|\tilde{S}_{s}f(g_{n_m}(y))-\tilde{S}_{s}f(y)\right|.\label{aaa4}
\end{eqnarray}
From (\ref{aaa}), (\ref{gndef}) and the continuity of $\tilde{S}_sf$, this expression converges to zero as $m\rightarrow \infty$. Hence we deduce that
\[\tilde{S}_{t+s}f(x)=\lim_{m\rightarrow\infty}{S}^{n_m}_{t}\tilde{S}_{s}f(x_{n_m})=\tilde{S}_t\tilde{S}_sf(x),\]
where the second equality holds because $Q^{n_m}(x_{n_m},t)\rightarrow \tilde{Q}(x,t)$ (by Lemma \ref{subseqlem}) and $\tilde{S}f\in C(M)$. This completes the proof of (\ref{semigroup}) for $f\in  \mathrm{Lip}_1(F)$ (again recalling that any such function can be extended to a function in $\mathrm{Lip}_1(M)$), and the extension to $f\in C(F)$ follows from an elementary approximation argument.

Finally, to prove the convergence of processes as at (\ref{felproc}), it will suffice to show convergence of finite-dimensional distributions, since we already have tightness from Proposition \ref{tightness}. This follows a standard argument (see \cite[Theorem 19.25]{Kall}, for example), though a little care is needed to account for the discrepancy between state spaces. Indeed, we are required to show that, for any $f_1,\dots,f_K\in C(M)$ and $0\leq t_1<\dots<t_K$,
\[\lim_{m\rightarrow \infty}E^{n_m}_{x_{n_m}}\left(\prod_{k\leq K}f_k(X^{n_m}_{t_k})\right)=\tilde{E}_{x}\left(\prod_{k\leq K}f_k(\tilde{X}_{t_k})\right).\]
We proceed by induction. For $K=1$, the result is clear from Lemma \ref{subseqlem}. Suppose the result holds for a product of $K-1$ terms. Applying the Markov property at $t_{K-1}$, we can write
\[E^{n_m}_{x_{n_m}}\left(\prod_{k\leq K}f_k(X^{n_m}_{t_k})\right)=
E^{n_m}_{x_{n_m}}\left(\prod_{k\leq K-1}f_k(X^{n_m}_{t_k})\times S^{n_m}_{t_K-t_{K-1}}f_K(X^{n_m}_{t_{K-1}})\right).\]
Observe that $\tilde{S}_{t_K-t_{K-1}}f_K\in C(F)$, and so can be extended to a function in $C(M)$. For this function, arguing as at (\ref{aaa4}), we have $\sup_{y\in F_{n_m}}|S^{n_m}_{t_K-t_{K-1}}f_K(y)-\tilde{S}_{t_K-t_{K-1}}f_K(y)|\rightarrow 0$. Hence
\begin{eqnarray*}
\lim_{m\rightarrow\infty}E^{n_m}_{x_{n_m}}\left(\prod_{k\leq K}f_k(X^{n_m}_{t_k})\right)&=&\lim_{m\rightarrow\infty}
E^{n_m}_{x_{n_m}}\left(\prod_{k\leq K-1}f_k(X^{n_m}_{t_k})\times \tilde{S}_{t_K-t_{K-1}}f_K(X^{n_m}_{t_{K-1}})\right)\\
&=&\tilde{E}_{x}\left(\prod_{k\leq K-1}f_k(\tilde{X}_{t_k})\times \tilde{S}_{t_K-t_{K-1}}f_K(\tilde{X}_{t_{K-1}})\right)\\
&=&\tilde{E}_{x}\left(\prod_{k\leq K}f_k(\tilde{X}_{t_k})\right),
\end{eqnarray*}
where the penultimate equality is an application of the inductive hypothesis.
\end{proof}

The main task in the remainder of the section is to show that any subsequential limit $\tilde{Q}$ is actually equal to $Q$. We do this by studying its resolvent, starting with the resolvent of a killed version of the process. In particular, for the associated Markov process, we define, for measurable $f:F\rightarrow\mathbb{R}_+$,
\[\tilde{G}_xf(y):=\tilde{E}_y\left(\int_0^{\tilde{\sigma}_x}f\left(\tilde{X}_s\right)ds\right),\qquad \forall x,y\in F,\]
where  $\tilde{\sigma}_x$ is the hitting time of $x$ by $\tilde{X}$, which may \emph{a priori} be infinite. In the following result, we show that $\tilde{G}_x$ is equal to $G_x$, as defined at (\ref{resolvdef}). The proof is similar to that of \cite[Proposition 5.6]{ALW}, but we avoid making an appeal to the `closed interval property' that was used there. We do, though, apply our resolvent approximation result from Lemma \ref{approxker}, and the resolvent convergence result of Proposition \ref{resolvprop}.

{\lem \label{killedresolv} Suppose Assumption \ref{mainassu} holds, and that $(F_n,R_n)$, ${n\geq 1}$, and $(F,R)$ are compact. If $\tilde{Q}$ is a subsequential limit of $(Q_n)_{n\geq 1}$, then
\[\tilde{G}_xf(y)=G_xf(y),\qquad \forall x,y\in F,\]
for every measurable $f:F\rightarrow\mathbb{R}_+$.}
\begin{proof} We start by showing $\tilde{G}_xf(y)\geq G_xf(y)$ for all $x,y\in F$ and positive $f\in C(F)$. Note that, because $F$ is closed in $M$, we can extend such a function to a positive function in $C(M)$. Suppose $(n_m)_{m\geq 1}$ is the subsequence giving rise to the limit $\tilde{Q}$. Let $x,y\in F$ and $(x_{n_m})_{m\geq 1}$, $(y_{n_m})_{m\geq 1}$ be sequences with $x_{n_m},y_{n_m}\in F_{n_m}$, $d_M(x_{n_m},x)\rightarrow 0$, and $d_M(y_{n_m},x)\rightarrow 0$. Then, applying that $X^{n_m}$ started from $y_{n_m}$ converges in distribution to $\tilde{X}$ started from $y$ (as established in Lemma \ref{fellerlem}), for any $\varepsilon>0$ we obtain
\begin{eqnarray*}
\tilde{G}_xf(y)&\geq &\limsup_{m\rightarrow \infty} {E}_{y_{n_m}}^{n_m}\left(\int_0^{{\sigma}^{n_m}_{\bar{B}_{n_m}(x_{n_m},\varepsilon)}}f\left({X}^{n_m}_s\right)ds\right)\\
&=&\limsup_{m\rightarrow \infty}G^{n_m}_{\bar{B}_{n_m}(x_{n_m},\varepsilon)}f(y_{n_m})\\
&=&\limsup_{m\rightarrow \infty}\int_{F_{n_m}}g^{n_m}_{\bar{B}_{n_m}(x_{n_m},\varepsilon)}(y_{n_m},z)f(z)\mu_{n_m}(dz)\\
&\geq &\limsup_{m\rightarrow \infty}\left(\int_{F_{n_m}}g^{n_m}_{x_{n_m}}(y_{n_m},z)f(z)\mu_{n_m}(dz)-2\varepsilon\|f\|_\infty\mu_{n_m}(F_{n_m})\right)\\
&=&\limsup_{m\rightarrow \infty}\left(G^{n_m}_{x_{n_m}}f(y_{n_m})-2\varepsilon\|f\|_\infty\mu_{n_m}(F_{n_m})\right)\\
&=&G_xf(y)-2\varepsilon\|f\|_\infty\mu(F),
\end{eqnarray*}
where for the second inequality we apply Lemma \ref{approxker}, and for the final equality we apply Proposition \ref{resolvprop}. Since $\varepsilon$ was arbitrary, this establishes the desired result.

Towards proving the reverse inequality, we start by showing that $\tilde{\sigma}_x=\tilde{\sigma}'_x$, $\tilde{P}_y$-a.s., where $\tilde{\sigma}'_x=\inf\{t>0:\:\tilde{X}_{t-}=x\}$. In particular, for any $\varepsilon,\delta,t_0>0$ we have that
\begin{equation}\label{c1}
\tilde{P}_y\left(R\left(x,\tilde{X}_{\tilde{\sigma}'_x}\right)\geq\varepsilon\right)
\leq\tilde{E}_y\left(\tilde{P}_{\tilde{X}_{\tilde{\sigma}_{\bar{B}_R(x,\delta)}}}\left(\sup_{s\leq t_0}R\left(x,\tilde{X}_{s}\right)\geq\varepsilon\right)\right)+\tilde{P}_y\left({\tilde{\sigma}'_x}-\tilde{\sigma}_{\bar{B}_R(x,\delta)}\geq t_0\right),
\end{equation}
where we have applied the fact that Feller processes are strong Markov (see \cite[Theorem 19.17]{Kall}, for example).  Now, by the right-continuity of $\tilde{X}$ at time zero (see \cite[Lemma 19.3]{Kall}), we have that
\[\tilde{P}_{x}\left(\sup_{s\leq t_0}R\left(x,\tilde{X}_{s}\right)>\varepsilon\right)\rightarrow 0,\]
as $t_0\rightarrow 0$. Hence, we can choose $t_0$ such that the probability on the left-hand side here is smaller that $\varepsilon/4$. It follows from this and \cite[Theorem 19.25]{Kall} that there exists a $\delta$ such that
\begin{equation}\label{c2}
\sup_{z\in \bar{B}_R(x,\delta)}\tilde{P}_{z}\left(\sup_{s\leq t_0}R\left(x,\tilde{X}_{s}\right)\geq \varepsilon\right)\leq \varepsilon/2.
\end{equation}
Moreover, it is clear from the definition that $\tilde{\sigma}_{\bar{B}_R(x,\delta)}\rightarrow {\tilde{\sigma}'_x}$ as $\delta\rightarrow 0$, $\tilde{P}_y$-a.s., and so, we can choose $\delta$ so that, in addition to \eqref{c2}, we also have
\begin{equation}\label{c3}
\tilde{P}_y\left({\tilde{\sigma}'_x}-\tilde{\sigma}_{\bar{B}_R(x,\delta)}\geq t_0\right)<\varepsilon/2.
\end{equation}
Applying \eqref{c2} and \eqref{c3} in \eqref{c1}, we obtain that
$\tilde{P}_y(R(x,\tilde{X}_{\tilde{\sigma}'_x})\geq\varepsilon)<\varepsilon$, and since $\varepsilon$ was arbitrary, this establishes the desired result. It thus follows that, in the setting of the first part of the proof
\begin{eqnarray*}
\tilde{G}_xf(y)&= &\tilde{E}_y\left(\int_0^{\tilde{\sigma}'_x}f\left(\tilde{X}_s\right)ds\right)\\
&= &\lim_{\varepsilon\rightarrow 0}\tilde{E}_y\left(\int_0^{\tilde{\sigma}_{\bar{B}_R(x,\varepsilon)}}f\left(\tilde{X}_s\right)ds\right)\\
&\leq&\lim_{\varepsilon\rightarrow 0}\liminf_{m\rightarrow \infty}{E}_{y_{n_m}}^{n_m}\left(\int_0^{{\sigma}^{n_m}_{\bar{B}_{n_m}(x_{n_m},\varepsilon)}}f\left({X}^{n_m}_s\right)ds\right)\\
&\leq&\liminf_{m\rightarrow \infty}G^{n_m}_{x_{n_m}}f(y_{n_m})\\
&=&G_xf(y),
\end{eqnarray*}
where the second inequality is clear from the definition of the resolvent (or Lemma \ref{approxker}), and the final equality is a consequence of Proposition \ref{resolvprop}.

We have so far shown that $\tilde{G}_xf(y)=G_xf(y)$ for all $x,y\in F$ and positive $f\in C(F)$. The extension to measurable $f:F\rightarrow\mathbb{R}_+$ is elementary.
\end{proof}

In the next result, we extend Lemma \ref{killedresolv} to hold for the $\alpha$-resolvent of the killed and unkilled processes. In particular, for $\alpha>0$, we write, for measurable $f:F\rightarrow\mathbb{R}_+$,
\[{G}^\alpha_xf(y):={E}_y\left(\int_0^{{\sigma}_x}e^{-\alpha s}f\left({X}_s\right)ds\right),\qquad \forall x,y\in F,\]
and also
\[{G}^\alpha f(y):={E}_y\left(\int_0^{\infty}e^{-\alpha s}f\left({X}_s\right)ds\right),\qquad \forall y\in F.\]
We define $\tilde{G}_x^\alpha$ and $\tilde{G}^\alpha$ from $\tilde{X}$ similarly.

{\lem\label{equalresolv} Suppose Assumption \ref{mainassu} holds, and that $(F_n,R_n)$, ${n\geq 1}$, and $(F,R)$ are compact. If $\tilde{Q}$ is a subsequential limit of $(Q_n)_{n\geq 1}$, then
\begin{equation}\label{r1}
\tilde{G}^\alpha_xf(y)=G^\alpha_xf(y),\qquad \forall x,y\in F,\:\alpha>0,
\end{equation}
for every measurable $f:F\rightarrow\mathbb{R}_+$. It moreover holds that
\begin{equation}\label{r2}
\tilde{G}^\alpha f(y)=G^\alpha f(y),\qquad \forall y\in F,\:\alpha>0,
\end{equation}
for every measurable $f:F\rightarrow\mathbb{R}_+$.}
\begin{proof} We start by noting that Lemma \ref{killedresolv} gives, for bounded, measurable $f:F\rightarrow\mathbb{R}_+$, and $x\in F$,
\[\|\tilde{G}_xf\|_{\infty}=\|{G}_xf\|_{\infty}\leq \|f\|_{\infty}\sup_{y\in F}E_y\sigma_x\leq \|f\|_{\infty} \mathrm{diam}_R(F)\mu(F),\]
where the final inequality follows from the commute time identity (Lemma \ref{commutetimeidentity}), with $\mathrm{diam}_R(F)$ the diameter of $(F,R)$, which is finite. Thus the Markov property for $\tilde{X}$ implies the following resolvent equation: for every $\alpha> 0$ and bounded, measurable $f:F\rightarrow \mathbb{R}_+$,
\[\tilde{G}^\alpha_xf(y)=\tilde{G}_xf(y)-\alpha\tilde{G}_x\tilde{G}^\alpha_xf(y),\qquad \forall x,y\in F,\]
with all the above terms being finite. Iterating this, we obtain that, for $\alpha<(\mathrm{diam}_R(F)\mu(F))^{-1}$,
\[\tilde{G}^\alpha_xf(y)=\sum_{i=1}^\infty (-\alpha)^{i-1}\tilde{G}^{\circ i}_xf(y),\qquad \forall x,y\in F,\]
where $\tilde{G}^{\circ i}_x$ is the $i$-fold composition of $\tilde{G}_x$. By the same argument, the corresponding result also holds for ${G}^\alpha_x$. Hence, Lemma \ref{killedresolv} yields that (\ref{r1}) holds for bounded, measurable $f:F\rightarrow \mathbb{R}_+$ and $\alpha<(\mathrm{diam}_R(F)\mu(F))^{-1}$. To extend to all $\alpha>0$, note that $\tilde{G}^\alpha_xf(y)$ is the Laplace transform of $\tilde{E}_y(f(\tilde{X}_t)\mathbf{1}_{t<\tilde{\sigma}_x})$ (considered as a function of $t$), and similarly ${G}^\alpha_xf(y)$ is the Laplace transform of ${E}_y(f({X}_t)\mathbf{1}_{t<{\sigma}_x})$. Since the two Laplace transforms agree on an open interval, they must agree for all $\alpha>0$. Extending to arbitrary measurable $f:F\rightarrow\mathbb{R}_+$ is elementary.

To establish (\ref{r2}), we apply the strong Markov property. In particular, fix $x_0,x_1\in F$ with $x_0\neq x_1$, and denote by ${\sigma}_0=\sigma_{x_0}$, and
\[\sigma_{i+1}:=\inf\left\{t>\sigma_i:\:X_t=x_0,\:x_1\in X_{[\sigma_i,t]}\right\}.\]
From the strong Markov property and the strong law of large numbers, we obtain that $\sigma_i\rightarrow \infty$, $P_y$-a.s. Hence, for bounded, measurable $f:F\rightarrow\mathbb{R}_+$, $\alpha>0$, $y\in F$,
\begin{eqnarray*}
G^\alpha f(y)&=&E_y\left(\int_0^{\sigma_0}e^{-\alpha s}f(X_s)ds\right)+\sum_{i=0}^\infty E_y\left(\int_{\sigma_i}^{\sigma_{i+1}}e^{-\alpha s}f(X_s)ds\right)\\
&=& G_{x_0}^\alpha f(y)+\sum_{i=0}^\infty
E_y\left(e^{-\alpha\sigma_0}\right)
E_{x_0}\left(e^{-\alpha\sigma_i}\right)E_{x_0}\left(\int_{0}^{\sigma_{1}}e^{-\alpha s}f(X_s)ds\right).
\end{eqnarray*}
Now, it is a simple computation that
\[E_y\left(e^{-\alpha\sigma_0}\right)=1-\alpha G^\alpha_{x_0}1(y),\]
\[E_{x_0}\left(e^{-\alpha\sigma_i}\right)=\left((1-\alpha G^\alpha_{x_1}1(x_0))(1-\alpha G^\alpha_{x_0}1(x_1))\right)^{i},\]
and
\[E_{x_0}\left(\int_{0}^{\sigma_{1}}e^{-\alpha s}f(X_s)ds\right)=G^\alpha_{x_1}f(x_0)+(1-\alpha G^\alpha_{x_1}1(x_0))G^\alpha_{x_0}f(x_1),\]
i.e.\ $G^\alpha f(y)$ can be expressed purely in terms of the killed resolvent. The same is true for the process $\tilde{X}$. Hence we obtain (\ref{r2}) from (\ref{r1}).
\end{proof}

Given what we have established so far, the proof of Proposition \ref{convtoQ} is almost immediate.

\begin{proof}[Proof of Proposition \ref{convtoQ}] Given Lemma \ref{subseqlem}, it will suffice to check that $\tilde{Q}=Q$, where $\tilde{Q}$ is any subsequential limit of $(Q_n)_{n\geq 1}$. Since both $(\tilde{S}_t)_{t\geq 0}$ and $(S_t)_{t\geq 0}$ are Feller (where the latter is the semigroup of $X$), we have from \cite[Theorem 19.4 and Lemma 19.5]{Kall} that they are uniquely determined by their resolvents $\tilde{G}^\alpha$ and ${G}^\alpha$ respectively (for any $\alpha>0$). Lemma \ref{equalresolv} gives us that the latter two operators are identical, and thus we have that so are $(\tilde{S}_t)_{t\geq 0}$ and $(S_t)_{t\geq 0}$, from which the result readily follows.
\end{proof}

\section{Proof of main result}\label{proofsec}

In view of the results of the previous two sections, it is now relatively straightforward to prove Theorem \ref{mainthm}. As previously, throughout this section, we suppose that Assumption \ref{mainassu} holds, and that $(F_n,R_n)$, ${n\geq 1}$, and $(F,R)$ are isometrically embedded into a common metric space $(M,d_M)$ in the way described in Lemma \ref{embeddings}. The next lemma gives the result for compact spaces.

{\lem\label{compactcase} If Assumption \ref{mainassu} holds, and $(F_n,R_n)$, ${n\geq 1}$, and $(F,R)$ are compact, then the conclusion of Theorem \ref{mainthm} holds.}
\begin{proof} This is an immediate consequence of Propositions \ref{tightness} and \ref{convtoQ}, cf.\ the subsequential convergence result established as part of Lemma \ref{fellerlem}.
\end{proof}

To extend to the locally compact case, we consider restrictions to bounded subsets. In particular, for $(F,R,\mu,\rho)\in\mathbb{F}$ and $r\geq 0$, introduce a continuous additive functional by setting
\[A^{(r)}_t:=\int_{F^{(r)}} L_t(x)\mu(dx),\]
and let $\tau^{(r)}(t):=\inf\{s:\:A^{(r)}_s>t\}$ be its right-continuous inverse. From Lemma \ref{restriction}, we have that
\begin{equation}\label{coupling}
\left(X^{(r)}_t\right)_{t\geq 0}:=\left(X_{\tau^{(r)}(t)}\right)_{t\geq 0}
\end{equation}
is the process associated with $(F^{(r)},R^{(r)},\mu^{(r)},\rho^{(r)})\in\mathbb{F}_c$. In the following lemma, we establish that $X^{(r)}$ converges in distribution to $X$ as $r\rightarrow\infty$.

{\lem \label{rconv} Let $(F,R,\mu,\rho)\in\mathbb{F}$. As $r\rightarrow\infty$,
\[P_\rho\left(\left(X^{(r)}_t\right)_{t\geq 0}\in \cdot\right)\rightarrow P_\rho\left( \left(X_t\right)_{t\geq 0}\in \cdot\right)\]
weakly as probability measures on $D(\mathbb{R}_+,F)$.}
\begin{proof} Let $\varepsilon,t>0$, and observe that
\[P_\rho\left(\sup_{s\leq t}R\left(X_s,X^{(r)}_s\right)>\varepsilon\right)\leq P_\rho\left(A_{t+1}^{(r)}\neq t+1\right)
\leq P_\rho\left(\sigma_{B_R(\rho,r)^c}<t+1\right).\]
Now, since $X$ is recurrent, we have from Lemma \ref{reclem} that $R(\rho,B_R(\rho,r)^c)\rightarrow\infty$. Hence we obtain from Lemma \ref{hittingtime}(b) that, for large $r$,
\begin{equation}\label{exit}
P_\rho\left(\sup_{s\leq t}R\left(X_s,X^{(r)}_s\right)>\varepsilon\right)\leq
4\left[
\frac{1}{R\left(\rho,{B}_R(\rho,r)^c\right)}+
\frac{t+1}{\mu(B_R(\rho,1))\left(R\left(\rho,{B}_R(\rho,r)^c\right)-1\right)}\right],
\end{equation}
where the right-hand side here converges to 0 as $r\rightarrow \infty$. The result follows.
\end{proof}

We next seek to extend the above result across elements of a sequence $(F_n,R_n,\mu_n,\rho_n)_{n\geq 1}$ in $\mathbb{F}$ satisfying Assumption \ref{mainassu}. We define $X^{n,r}$ from $X^n$ in the way $X^{(r)}$ was defined from $X$ at (\ref{coupling}), and note this is the process naturally associated with $(F_n^{(r)},R_n^{(r)},\mu_n^{(r)},\rho_n^{(r)})$.

{\lem\label{nrn} Suppose Assumption \ref{mainassu} holds. It is then the case that, for any $\varepsilon, t>0$,
\[\lim_{r\rightarrow \infty}\limsup_{n\rightarrow\infty}P^n_{\rho_n}\left(\sup_{s\leq t}R_n\left(X^n_s,X^{n,r}_{s}\right)>\varepsilon\right)=0.\]}
\begin{proof} As at (\ref{exit}), we have that
\[P^n_{\rho_n}\left(\sup_{s\leq t}R_n\left(X^n_s,X^{n,r}_s\right)>\varepsilon\right)\leq
4\left[
\frac{1}{R_n\left(\rho,{B}_n(\rho_n,r)^c\right)}+
\frac{t+1}{\mu_n(B_n(\rho_n,1))\left(R_n\left(\rho,{B}_n(\rho_n,r)^c\right)-1\right)}\right]\]
whenever $R_n(\rho_n,B_n(\rho_n,r)^c)>1$. Since $\liminf_{n\rightarrow\infty}\mu_n(B_n(\rho_n,1))\geq \mu(B_R(\rho,1))>0$, the lemma thus follows from the resistance growth condition at \eqref{resgrowthcond}.
\end{proof}

Proceeding exactly as for Lemma \ref{compactcase}, we can also obtain the convergence of the restricted processes. The only additional point to note in establishing the following lemma is that the embeddings of Lemma \ref{embeddings} can be chosen in such a way that the result holds simultaneously for Lebesgue almost-every $r\geq 0$.

{\lem \label{compactconv} Suppose Assumption \ref{mainassu} holds. It is then the case that, for Lebesgue almost-every $r\geq 0$,
\[P^n_{\rho_n}\left(\left(X^{n,r}_t\right)_{t\geq 0}\in \cdot\right)\rightarrow P_\rho\left( \left(X^{(r)}_t\right)_{t\geq 0}\in \cdot\right)\]
weakly as probability measures on $D(\mathbb{R}_+,M)$.}
\bigskip

From Lemmas \ref{rconv}, \ref{nrn} and \ref{compactconv}, the conclusion of Theorem \ref{mainthm} is standard (see \cite[Theorem 4.28]{Kall}, for example).

\section{Spatial embeddings}\label{embsec}

As raised in Remark \ref{explosionrem}(c), it might sometimes be natural to consider the convergence of processes with respect to a different topology than that generated by the resistance metric. In particular, the spaces of interest might already be embedded into some common space, upon which one might hope to prove a convergence statement. In this section, we present a result aimed at tackling this problem, in the case when the relevant embeddings are continuous with respect to the resistance metric (see Theorem \ref{embthm} below), and also a variation of this for random spaces (Theorem \ref{randembthm}). In particular, we now consider the space of spatial, measured resistance metric spaces $\mathbb{F}^*$, which is defined to be the collection of quintuplets $(F,R,\mu,\rho,\phi)$ such that $(F,R,\mu,\rho)\in\mathbb{F}$, and $\phi:F\rightarrow(N,d_N)$ is a continuous map, where the image space $(N,d_N)$ is a complete, separable metric space (fixed across all elements of $\mathbb{F}^*$). Note that the terminology `spatial' is borrowed from \cite{DLeG}, where a similar space of real trees with an embedding into Euclidean space was introduced.

To define a topology on $\mathbb{F}^*$, it will be convenient to start by considering the larger space $\mathbb{K}^*$ of elements $(K,d_K,\mu,\rho,\phi)$, where now we simply assume: $(K,d_K)$ is a metric space in which closed bounded sets are compact, $\mu$ is a locally finite Borel measure on $K$, $\rho$ is a marked point of $K$, and $\phi$ is a continuous map from $K$ into $N$. Similarly to Section \ref{ghvsubsec}, we first consider the case when $(K,d_K)$ is compact; write $\mathbb{K}^*_c$ for the subset of $\mathbb{K}^*$ for which this is the case. For two elements of $\mathbb{K}^*_c$, we set $\Delta^*_c((K,d_K,\mu,\rho,\phi),(K',d_{K'},\mu',\rho',\phi'))$ to be equal to
\begin{equation*}
\inf_{\substack{M,\psi,\psi',\mathcal{C}:\\(\rho,\rho')\in\mathcal{C}}}
\left\{d_M^P\left(\mu\circ\psi^{-1},\mu'\circ\psi'^{-1}\right)+
\sup_{(x,x')\in\mathcal{C}}\left(d_M\left(\psi(x),\psi'(x')\right)+d_N\left(\phi(x),\phi'(x')\right)\right)\right\},
\end{equation*}
where the infimum is taken over all metric spaces $(M,d_M)$, isometric embeddings $\psi:(K,d_K)\rightarrow M$, $\psi':(K',d_{K'})\rightarrow M$, and correspondences $\mathcal{C}$ between $K$ and $K'$. Note that, by a correspondence $\mathcal{C}$ between $K$ and $K'$, we mean a subset of $K\times K'$ such that for every $x\in K$ there exists at least one $x'\in K'$ such that $(x,x')\in\mathcal{C}$ and conversely for every $x'\in K'$ there exists at least one $x\in K$ such that $(x,x')\in \mathcal{C}$. We note that it is possible to check that $(\mathbb{K}^*_c,\Delta^*_c)$ is a separable metric space exactly as in \cite[Proposition 3.1]{BCK} (cf.\ \cite[Lemma 2.1]{CHKmix}). More generally, for elements of $\mathbb{K}^*$, we consider restrictions to compact subsets. Specifically, given $(K,d_K,\mu,\rho,\phi)\in\mathbb{K}^*$, write $K^{(r)}$ for the closure of the ball of radius $r$ about $\rho$, i.e.\ $K^{(r)}:=\bar{B}_K(\rho,r)$, $d_K^{(r)}:=d_K|_{K^{(r)}\times K^{(r)}}$, $\mu^{(r)}:=\mu(\cdot\cap K^{(r)})$, $\rho^{(r)}:=\rho$, and $\phi^{(r)}:=\phi|_{K^{(r)}}$. We then say $(K_n,d_{K_n},\mu_n,\rho_n,\phi_n)\rightarrow (K,d_K,\mu,\rho,\phi)$ in the spatial Gromov-Hausdorff-vague topology if and only if
\[\Delta^*_c\left((K_n^{(r)},d^{(r)}_{K_n},\mu^{(r)}_n,\rho^{(r)}_n,\phi^{(r)}_n),(K^{(r)},d^{(r)}_{K},\mu^{(r)},\rho^{(r)},\phi^{(r)})\right)\rightarrow 0\]
for Lebesgue almost-every $r\geq 0$.

We are now in a position to state our main conclusion for stochastic processes associated with spatial, measured resistance metric spaces.

{\thm \label{embthm} If a sequence $(F_n,R_n,\mu_n,\rho_n,\phi_n)_{n\geq 1}$ in $\mathbb{F}^*$ satisfies
\[\left(F_n,R_n,\mu_n,\rho_n,\phi_n\right)\rightarrow\left(F,R,\mu,\rho,\phi\right)\]
in the spatial Gromov-Hausdorff-vague topology for some $(F,R,\mu,\rho,\phi)\in\mathbb{F}^*$, and (\ref{resgrowthcond}) holds, then
\[P^n_{\rho_n}\left(\left(\phi_n\left(X^n_t\right)\right)_{t\geq 0}\in \cdot\right)\rightarrow P_{\rho}\left(\left(\phi\left(X_t\right)\right)_{t\geq 0}\in \cdot\right)\]
weakly as probability measures on $D(\mathbb{R}_+,N)$.}
\bigskip

Before proving this result, let us also state a version that is useful for random spaces. In particular, in the following result, we suppose that we have random elements $(F_n,R_n,\mu_n,\rho_n,\phi_n)_{n\geq 1}$ and $(F,R,\mu,\rho,\phi)$ of $\mathbb{F}^*$, built on a probability space with probability measure $\mathbf{P}$. We will show that if $(F_n,R_n,\mu_n,\rho_n,\phi_n)\rightarrow(F,R,\mu,\rho,\phi)$ in a distributional sense with respect to the spatial Gromov-Hausdorff-vague topology, and a probabilistic version of (\ref{resgrowthcond}) holds, then we also have convergence of embedded stochastic processes under the annealed measure, that is the probability measure on $D(\mathbb{R}_+,N)$ obtained by integrating out the randomness of the state spaces, i.e. the annealed law of the limit process is given by
\[\mathbb{P}_{\rho}\left(\left(\phi\left(X_t\right)\right)_{t\geq 0}\in \cdot\right):=
\int{P}_{\rho}\left(\left(\phi\left(X_t\right)\right)_{t\geq 0}\in \cdot\right)d\mathbf{P},\]
and we define $\mathbb{P}_{\rho_n}^n$, $n\geq 1$, similarly.

{\thm \label{randembthm} Suppose $(F_n,R_n,\mu_n,\rho_n,\phi_n)_{n\geq 1}$ and $(F,R,\mu,\rho,\phi)$ are random elements of $\mathbb{F}^*$ such that
\begin{equation}\label{restrand}
\left(F_n^{(r)},R_n^{(r)},\mu_n^{(r)},\rho_n^{(r)},\phi_n^{(r)}\right)\buildrel{d}\over\rightarrow\left(F^{(r)},R^{(r)},\mu^{(r)},\rho^{(r)},\phi^{(r)}\right)
\end{equation}
in the spatial Gromov-Hausdorff-vague topology for Lebesgue almost-every $r\geq 0$, and
\begin{equation}\label{randresgrowthcond}
\lim_{r\rightarrow\infty}\liminf_{n\rightarrow\infty}\mathbf{P}\left(R_n\left(\rho_n,B_n\left(\rho_n,r\right)^c\right)\geq \lambda\right)=1,\qquad \forall\lambda\geq 0,
\end{equation}
then
\[\mathbb{P}^n_{\rho_n}\left(\left(\phi_n\left(X^n_t\right)\right)_{t\geq 0}\in \cdot\right)\rightarrow \mathbb{P}_{\rho}\left(\left(\phi\left(X_t\right)\right)_{t\geq 0}\in \cdot\right)\]
weakly as probability measures on $D(\mathbb{R}_+,N)$.}
\bigskip

The following lemma provides the key technical points that will be needed for the proof of Theorem \ref{embthm}. Firstly, it incorporates correspondences and spatial embeddings to Lemma \ref{embeddings}. Secondly, it gives pointwise convergence and equicontinuity for spatial embeddings.

{\lem (a) Suppose $(F_n,R_n,\mu_n,\rho_n,\phi_n)$, $n\geq 1$, and $(F,R,\mu,\rho,\phi)$ are elements of $\mathbb{F}^*_c$ such that $(F_n,R_n,\mu_n,\rho_n,\phi_n)\rightarrow(F,R,\mu,\rho,\phi)$ in the spatial Gromov-Hausdorff-vague topology. It is then possible to embed $(F_n,R_n)$, $n\geq 1$, and $(F,R)$ isometrically into the same compact metric space $(M,d_M)$ in such a way that, for Lebesgue almost-every $r\geq 0$, the claims at (\ref{embeddingseq}) hold. Moreover, it is simultaneously possible to find correspondences $\mathcal{C}_n\subseteq F\times F_n$ such that
\begin{equation}\label{v1}
\varepsilon_n^{(1)}:=\sup_{(x,y)\in\mathcal{C}_n}d_M(x,y)\rightarrow 0,
\end{equation}
\begin{equation}\label{v2}
\varepsilon_n^{(2)}:=\sup_{(x,y)\in\mathcal{C}_n}d_N(\phi(x),\phi_n(y))\rightarrow 0,
\end{equation}
where we have identified the various objects with their embeddings into $M$.\\
(b) In the setting of part (a), we also have that if $x_n\in F_n$, $n\geq 1$, and $x\in F$ are such that $d_M(x_n,x)\rightarrow0$, then
\begin{equation}\label{pointwise}
d_N(\phi_n(x_n),\phi(x))\rightarrow 0.
\end{equation}
Moreover,
\begin{equation}\label{phiequi}
\lim_{\delta\rightarrow0}\limsup_{n\rightarrow\infty}\sup_{\substack{x,y\in F_n:\\R_n(x,y)\leq \delta}}d_N(\phi_n(x),\phi_n(y))=0.
\end{equation}}
\begin{proof} An embedding satisfying the requirements of part (a) of the lemma is described (in a very slightly modified setting, with no marked point, and where the functions $\phi$ are replaced by a function with two arguments plus a time parameter) in the proof of \cite[Lemma 2.2]{CHKmix}. We now check part (b). First, suppose $x_n\in F_n$, $n\geq 1$, and $x\in F$ are such that $d_M(x_n,x)\rightarrow0$. Let $x_n'\in F$ be such that $(x_n',x_n)\in \mathcal{C}_n$. Then
\begin{eqnarray*}
d_N(\phi_n(x_n),\phi(x))&\leq& d_N(\phi_n(x_n),\phi(x_n'))+d_N(\phi(x_n'),\phi(x))\\
&\leq & \varepsilon_n^{(2)}+\sup_{\substack{y,z\in F:\\ R(y,z)\leq d_M(x,x_n)+\varepsilon_n^{(1)}}}d_N(\phi(y),\phi(z)),
\end{eqnarray*}
which converges to zero by (\ref{v1}), (\ref{v2}), and the continuity of $\phi$. This establishes (\ref{pointwise}). For (\ref{phiequi}), we proceed similarly. In particular, given $x,y\in F_n$ with $R_n(x,y)\leq \delta$, choose $x',y'\in F$ such that $(x',x),(y',y)\in \mathcal{C}_n$, then
\begin{eqnarray*}
d_N(\phi_n(x),\phi_n(y))&\leq &d_N(\phi_n(x),\phi(x'))+d_N(\phi(x'),\phi(y'))+d_N(\phi(y'),\phi_n(y))\\
&\leq &2\varepsilon_n^{(2)}+\sup_{\substack{y,z\in F:\\ R(y,z)\leq \delta+2\varepsilon_n^{(1)}}}d_N(\phi(y),\phi(z)),
\end{eqnarray*}
from which (\ref{phiequi}) follows by again applying (\ref{v1}), (\ref{v2}), and the continuity of $\phi$.
\end{proof}

With these preparations in place, we are ready to prove the main results of this section.

\begin{proof}[Proof of Theorem \ref{embthm}] To begin with, suppose $(F_n,R_n,\mu_n,\rho_n,\phi_n)$, $n\geq 1$, and $(F,R,\mu,\rho,\phi)$ are elements of $\mathbb{F}^*_c$. From the argument used to prove Theorem \ref{mainassu}, we know that
\[P^n_{\rho_n}\left(\left(X^n_t\right)_{t\geq 0}\in \cdot\right)\rightarrow P_{\rho}\left(\left(X_t\right)_{t\geq 0}\in \cdot\right)\]
weakly as probability measures on $D(\mathbb{R}_+,M)$. Consequently, the finite-dimensional distributions converge, meaning that, for any $0\leq t_1<\dots<t_K$, we have that $(X^n_{t_k})_{k=1}^K\rightarrow (X_{t_k})_{k=1}^K$ in distribution in $M^K$. From (\ref{pointwise}), it thus follows that
$(\phi_n(X^n_{t_k}))_{k=1}^K\rightarrow (\phi(X_{t_k}))_{k=1}^K$ in distribution in $N^K$. To extend from this finite-dimensional distribution statement to convergence in $D(\mathbb{R}_+,N)$, we need to check tightness. For this, it will suffice to check the analogue of Lemma \ref{fluctuni}. To this end, observe that, under $P_x^n$,
\[\sup_{s\leq t}d_N(\phi_n(x),\phi_n(X^n_s))\leq \sup_{\substack{y\in F_n:\\R_n(x,y)\leq \sup_{s\leq t}R_n(x,X^n_s)}}d_N(\phi_n(x),\phi_n(y)).\]
Hence from Lemma \ref{fluctuni} and \eqref{phiequi}, we obtain that
\[\lim_{t\rightarrow0}\limsup_{n\rightarrow \infty}\sup_{x\in F_n}P^n_x\left( \sup_{s\leq t}d_N(\phi_n(x),\phi_n(X^n_s))\geq\varepsilon\right)=0,\]
as desired. In particular, we have now shown in the compact case that
\[P^n_{\rho_n}\left(\left(\phi_n\left(X^n_t\right)\right)_{t\geq 0}\in \cdot\right)\rightarrow P_{\rho}\left(\left(\phi\left(X_t\right)\right)_{t\geq 0}\in \cdot\right)\]
weakly as probability measures on $D(\mathbb{R}_+,N)$. One can extend this result to the locally compact case by exactly the same strategy as was applied in the proof of Theorem \ref{mainthm}, namely the consideration of the traces of the relevant processes on compact subsets.
\end{proof}

\begin{proof}[Proof of Theorem \ref{randembthm}] First, suppose $r>0$ is chosen such that (\ref{restrand}) holds. Since $(\mathbb{K}_c^*,\Delta_c^*)$ is separable, it is possible to find realisations of the spaces $(F_n^{(r)},R_n^{(r)},\mu_n^{(r)},\rho_n^{(r)},\phi_n^{(r)})$, $n\geq 1$, and $(F^{(r)},R^{(r)},\mu^{(r)},\rho^{(r)},\phi^{(r)})$ so that the convergence takes place almost-surely \cite[Theorem 4.30]{Kall}. For these realisations, we obtain from Theorem \ref{embthm} that, almost-surely,
\[P^n_{\rho_n}\left(\left(\phi_n\left(X^{n,r}_t\right)\right)_{t\geq 0}\in \cdot\right)\rightarrow P_{\rho}\left(\left(\phi\left(X^{(r)}_t\right)\right)_{t\geq 0}\in \cdot\right)\]
weakly as probability measures on $D(\mathbb{R}_+,N)$, where $X^{n,r}$ and $X^{(r)}$ are the processes introduced in Section \ref{proofsec}. Integrating out over the state spaces thus yields
\begin{equation}\label{yo1}
\mathbb{P}^n_{\rho_n}\left(\left(\phi_n\left(X^{n,r}_t\right)\right)_{t\geq 0}\in \cdot\right)\rightarrow \mathbb{P}_{\rho}\left(\left(\phi\left(X^{(r)}_t\right)\right)_{t\geq 0}\in \cdot\right)
\end{equation}
weakly as probability measures on $D(\mathbb{R}_+,N)$. Next, similarly to the proof of Lemma \ref{nrn}, we have that, for any $r_0>0$,
\begin{eqnarray*}
\lefteqn{\mathbb{P}^n_{\rho_n}\left(\sup_{s\leq t}d_N\left(\phi_n(X^n_s),\phi_n(X^{n,r}_s)\right)>\varepsilon\right)}\\
&\leq&
\mathbf{E}\left[\min\left\{1,4\left(
\frac{2r_0}{\lambda}+
\frac{t+1}{\mu_n(B_n(\rho_n,2r_0))\left(\lambda-2r_0\right)}\right)\right\}\right]+\mathbf{P}\left(R_n\left(\rho_n,B_n\left(\rho_n,r\right)^c\right)<\lambda\right).
\end{eqnarray*}
In particular, choosing $r_0$ so that (\ref{restrand}) holds, we have that $\mu_n(\bar{B}_n(\rho_n,r_0))\rightarrow \mu(\bar{B}_R(\rho,r_0))$, and so we obtain
\begin{eqnarray*}
\lefteqn{\limsup_{n\rightarrow\infty}\mathbb{P}^n_{\rho_n}\left(\sup_{s\leq t}d_N\left(\phi_n(X^n_s),\phi_n(X^{n,r}_s)\right)>\varepsilon\right)}\\
&\leq&
\mathbf{E}\left[\min\left\{1,4\left(
\frac{2r_0}{\lambda}+
\frac{t+1}{\mu(\bar{B}_R(\rho,r_0))\left(\lambda-2r_0\right)}\right)\right\}\right]+\limsup_{n\rightarrow\infty}\mathbf{P}\left(R_n\left(\rho_n,B_n\left(\rho_n,r\right)^c\right)< \lambda\right).
\end{eqnarray*}
Letting $r\rightarrow\infty$ and then $\lambda\rightarrow\infty$ (and applying (\ref{randresgrowthcond})), we thus deduce
\begin{equation}\label{yo2}
\lim_{r\rightarrow\infty}\limsup_{n\rightarrow\infty}\mathbb{P}^n_{\rho_n}\left(\sup_{s\leq t}d_N\left(\phi_n(X^n_s),\phi_n(X^{n,r}_s)\right)>\varepsilon\right)=0.
\end{equation}
The same argument yields that
\begin{equation}\label{yo3}
\lim_{r\rightarrow\infty}\mathbb{P}_{\rho}\left(\sup_{s\leq t}d_N\left(\phi(X_s),\phi(X^{(r)}_s)\right)>\varepsilon\right)=0.
\end{equation}
Putting together (\ref{yo1}), (\ref{yo2}) and (\ref{yo3}) completes the proof.
\end{proof}

\section{Examples}\label{exsec}

\subsection{Trees}

Since the shortest path metric on a tree is a resistance metric, such objects fit naturally into the framework of this article. With the scaling of stochastic processes on trees now being well-understood, however, we will not present details of particular cases here. Rather, to illustrate the applicability of the results, we simply list some previous work on scaling stochastic processes on trees that the conclusion of \cite{ALW} and this article covers (see the cited references for further details).
\begin{itemize}
  \item Critical Galton-Watson trees, conditioned on their size, both finite variance \cite{Crofin} and infinite variance \cite{Croyrogt}. Also, the incipient infinite cluster of the critical Galton-Watson tree with finite variance \cite[Section 7.4]{ALW} (see \cite{Kesten} for the initial work on the random walk on this object). In \cite{ALW}, the discussion further covers the so-called `height process' of the latter random walk model, that is, the process which measures the distance from the root of the tree to the walker; by taking spatial embeddings (into $\mathbb{R}_+$) given by the height function, one can obtain convergence of height processes as an application of Theorem \ref{randembthm}. One might alternatively take as the spatial embedding the height of the `projection to the backbone' of the random walk (note the incipient infinite cluster of a critical Galton-Watson tree almost-surely has a unique path to infinity, which is referred to as its backbone), and apply Theorem \ref{randembthm} to recover the scaling result proved in \cite{bACF0} for this function of the random walk (in the aforementioned article, a detailed characterisation of the limiting process is given in terms of a so-called `spatially-subordinated Brownian motion').
  \item The (non-lattice) branching random walk, where the underlying tree is a critical Galton-Watson tree with exponential tails for the offspring distribution, and the steps have a centred, continuous distribution with fourth order polynomial tail decay \cite{CroyHd}. (This is an example to which Theorem \ref{randembthm} applies, with embeddings into $\mathbb{R}^d$.)
  \item (Subsequentially,) the uniform spanning tree in two dimensions \cite{BCK}. (With its canonical embedding into $\mathbb{R}^2$, this is another example to which Theorem \ref{randembthm} applies.)
  \item $\Lambda$-coalescent measure trees \cite[Section 7.5]{ALW}.
\end{itemize}
Of course, there are also many other interesting examples of families of trees known to converge under scaling with respect to the Gromov-Hausdorff-vague topology for which convergence of processes would also follow from \cite{ALW} or Theorem \ref{mainthm}. These include the model of Markov branching trees \cite{HM}, and the minimal spanning tree of the complete graph \cite{ABGM}.

A potential advantage of the wider setting of this article compared to previous work on trees is that it gives a route for proving similar results for graphs that are not trees, but which have a tree as a scaling limit. Indeed, one model for which an appealing conjecture can be made is the incipient infinite cluster of bond percolation on $\mathbb{Z}^d$ in high dimensions, that is when $d>6$. In particular, it is expected that this model satisfies the same scaling properties as branching random walk, and thus one might anticipate that if $\mathrm{IIC}$ is the incipient infinite cluster (see \cite{HJ} for a construction), $R_{\mathrm{IIC}}$ is the resistance metric on this (when individual edges have unit resistance), $\mu_{\mathrm{IIC}}$ is the counting measure on $\mathrm{IIC}$ (placing mass one on each vertex), and $I_{\mathrm{IIC}}:\mathrm{IIC}\rightarrow\mathbb{R}^d$ is the identity map, then $(\mathrm{IIC},R_{\mathrm{IIC}},\mu_{\mathrm{IIC}},0,I_{\mathrm{IIC}})\in\mathbb{F}^*$, and the rescaled sequence
\[\left(\mathrm{IIC},n^{-2}R_{\mathrm{IIC}},n^{-4}\mu_{\mathrm{IIC}},0,n^{-1}I_{\mathrm{IIC}}\right),\qquad n\geq 1,\]
satisfies the conditions of Theorem \ref{randembthm}, with limit being (an unbounded version of) the continuum random tree mapped into Euclidean space by a tree-indexed Brownian motion (or, to state this another way, the integrated super-Brownian excursion), cf.\ the conjecture of \cite{CroyHd}, where the Brownian motion on the limiting space was constructed for $d\geq 8$. Note that early work connecting the incipient infinite cluster and the integrated super-Brownian excursion was accomplished in \cite{HS}. See also \cite{HHH} for up-to-date work on random walks and resistance on the incipient infinite cluster in high dimensions, which supports the above scaling picture. Finally, we note that, in this direction, the corresponding limit for the random walk on the lattice branching random walk has recently been established in \cite{bACF2}. The latter work depends on a general scaling theorem proved in \cite{bACF1}, which, similarly to the present work, depends critically on the connections between the resistance metric and stochastic processes.

\subsection{Range of random walk}

Continuing the discussion from the previous section, we now present a simple example of a random graph that is asymptotically a tree, and for which we can check the conditions of Theorems \ref{embthm} and \ref{randembthm} hold; this is the random walk on the range of a random walk in high dimensions. More precisely, for $(S_n)_{n\geq 0}$ a simple random walk on $\mathbb{Z}^d$ started from 0 (built on a probability space with probability measure $\mathbf{P}$), set $\mathcal{G}=(V(\mathcal{G}),E(\mathcal{G}))$ to be the graph with vertex set
\[V\left(\mathcal{G}\right):=\left\{S_n:\:n\geq 0\right\},\]
and edge set
\[E\left(\mathcal{G}\right):=\left\{\{S_n,S_{n+1}\}:\:n\geq 0\right\}.\]
For $d\geq5$, it is known that this graph is essentially one-dimensional, and that the associated random walk behaves like simple random walk on a half-line \cite{RWRRW}. Here we explain how to reprove the main results of the latter article using our resistance form approach (since the result is already known, we attempt to keep the exposition brief, and so omit detailed computations). We will denote by: $R_\mathcal{G}$ the resistance on $\mathcal{G}$, supposing individual edges have unit resistance; $\mu_\mathcal{G}$ for the measure on $V(\mathcal{G})$ placing mass one on each vertex; and $d_\mathcal{G}$ the shortest path graph distance on $\mathcal{G}$. In particular, with $I_\mathcal{G}$ being the identity map from $\mathcal{G}$ to $\mathbb{R}^d$, both  $(\mathcal{G},R_\mathcal{G},\mu_\mathcal{G},0,d_\mathcal{G}(0,\cdot))$
and
$(\mathcal{G},R_\mathcal{G},\mu_\mathcal{G},0,I_\mathcal{G})$
are ($\mathbf{P}$-a.s.) in versions of $\mathbb{F}^*$ (with embeddings into $\mathbb{R}_+$ and $\mathbb{R}^d$ respectively). NB.\ The trace of a random walk is always a recurrent graph \cite{Trace}.

To capture the one-dimensional nature of the graph $\mathcal{G}$, it is useful to introduce the collection of cut-times of $S$, that is the set
\[\mathcal{T}:=\left\{n:\:S_{[0,n]}\cap S_{[n+1,\infty)}\right\}.\]
One has that, $\mathbf{P}$-a.s., $\mathcal{T}$ is infinite, and so we can write its elements as $0\leq T_1<T_2<\dots$. The associated set of cut-points is given by $\mathcal{C}:=\{n\geq 1:\:C_n\}$, where $C_n:=S_{T_n}$. Applying laws of large numbers for various graph quantities measured at cut-times, we are able to deduce the following. Note that we write: $d_E$ for the Euclidean metric on $\mathbb{R}_+$; $\mathcal{L}$ for Lebesgue measure on $\mathbb{R}_+$; and $(W_t)_{t\geq 0}$ for a $d$-dimensional Brownian motion, started from $0$, also built on the probability space with probability measure $\mathbf{P}$.

{\propn For some deterministic constants $c_1,c_2,c_3,c_4\in(0,\infty)$, the following statements hold.\\
(a) $\mathbf{P}$-a.s., the sequence
\[\left(\mathcal{G},n^{-1}R_\mathcal{G},n^{-1}\mu_\mathcal{G},0,n^{-1}d_\mathcal{G}(0,\cdot)\right),\qquad n\geq 1,\]
satisfies the conditions of Theorem \ref{embthm}, with limit given by $(\mathbb{R}_+,c_1d_E,c_2\mathcal{L},0,c_3d_E(0,\cdot))$.\\
(b) Under $\mathbf{P}$, the sequence
\[\left(\mathcal{G},n^{-1}R_\mathcal{G},n^{-1}\mu_\mathcal{G},0,n^{-1/2}I_\mathcal{G}\right),\qquad n\geq 1,\]
satisfies the conditions of Theorem \ref{randembthm}, with (random) limit given by $(\mathbb{R}_+,c_1d_E,c_2\mathcal{L},0, c_4W_{\cdot})$.}
\begin{proof} We first claim that there exist deterministic constants such that, $\mathbf{P}$-a.s.,
\begin{equation}\label{ergodiclims}
\frac{T_n}{n}\rightarrow c_1,\qquad \frac{R_\mathcal{G}(0,C_n)}{n}\rightarrow c_2,\qquad \frac{\mu_\mathcal{G}(\{S_0,S_1,\dots,S_{T_n}\})}{n}\rightarrow c_3,\qquad \frac{d_\mathcal{G}(0,C_n)}{n}\rightarrow c_4.
\end{equation}
To prove this, we observe that the results corresponding to the first, second and fourth of these limits for the two-sided simple random walk $(\tilde{S}_n)_{n\in \mathbb{Z}}$ (that is, where $(\tilde{S}_n)_{n\geq 0}$ and $(\tilde{S}_{-n})_{n\geq 0}$ are independent simple random walks started from 0) readily follows from the ergodic theory argument of \cite[Lemmas 2.1 and 2.2]{RWRRW}. Moreover, the result corresponding to the third limit in the two-sided case can be established in an identical fashion. Note that, in the two-sided case, the set of cut-times is given by $\tilde{\mathcal{T}}:=\{n:\:\tilde{S}_{(-\infty,n]}\cap \tilde{S}_{[n+1,\infty)}\}$, the elements of which can be written as $\dots<\tilde{T}_{-1}<\tilde{T}_0<0\leq \tilde{T}_1<\tilde{T}_2<\dots$, and the set of cut-points $\tilde{\mathcal{C}}:=\{n\in\mathbb{Z}:\:\tilde{C}_n\}$, where $\tilde{C}_n:=\tilde{S}_{\tilde{T}_n}$. Identifying $(S_n)_{n\geq 0}$ and $(\tilde{S}_n)_{n\geq 0}$, one clearly has that $\mathbf{P}$-a.s., $C_n:=\tilde{C}_{n-n_0}$ for some (random, but almost-surely finite) $n_0$. Applying this, it is elementary to deduce the limits for the one-sided process, as stated at (\ref{ergodiclims}). Moreover, given (\ref{ergodiclims}), it is a simple exercise to verify part (a) of the lemma. To prove part (b) of the lemma, one only needs to incorporate the fact that $(n^{-1/2}S_{nt})_{t\geq0}\rightarrow (W_t)_{t\geq 0}$ in distribution under $\mathbf{P}$.
\end{proof}

In view of this proposition, Theorems \ref{embthm} and \ref{randembthm} immediately give convergence of the corresponding processes with respect to the quenched and annealed laws respectively. We write $((X^\mathcal{G}_t)_{t\geq 0},(P_x^\mathcal{G})_{x\in V(\mathcal{G})})$ for the process associated with $(\mathcal{G},R_\mathcal{G},\mu_\mathcal{G},0)\in\mathbb{F}$. Moreover, we denote by $\mathbb{P}_0^\mathcal{G}:=\int P_0^\mathcal{G}d\mathbf{P}$ the corresponding annealed law (when the process is viewed as a random element of $D(\mathbb{R}_+,\mathbb{R}^d)$). To describe the limiting process, we let $(B_t)_{t\geq 0}$ be a one-dimensional Brownian motion, so that $(|B_t|)_{t\geq 0}$ is the process associated with $(\mathbb{R}_+,d_E,\mathcal{L},0)\in\mathbb{F}$, and $P_0$ be its law. Finally, we set $\mathbb{P}_0:=\mathbf{P}\times P_0$.

{\cor[{cf.\ \cite[Theorem 1.1]{RWRRW}}] For some deterministic constants $c_1,c_2\in(0,\infty)$, the following statements hold.\\
(a) $\mathbf{P}$-a.s.,
\[P^\mathcal{G}_0\left(\left(n^{-1}d_\mathcal{G}(0,X^\mathcal{G}_{tn^2})\right)_{t\geq 0}\in \cdot\right)\rightarrow P_0\left(\left(|B_{c_1t}|\right)_{t\geq 0}\in \cdot\right)\]
weakly as probability measures on $D(\mathbb{R}_+,\mathbb{R}_+)$.\\
(b)
\[\mathbb{P}_0^\mathcal{G}\left(\left(n^{-1}X^\mathcal{G}_{tn^2}\right)_{t\geq 0}\in \cdot\right)\rightarrow \mathbb{P}_0\left(\left(W_{|B_{c_2t}|}\right)_{t\geq 0}\in \cdot\right)\]
weakly as probability measures on $D(\mathbb{R}_+,\mathbb{R}^d)$.}

\subsection{Fused resistance forms}

Another example of a family of random graphs that are nearly trees, and which the results of this article are applicable is the Erd\H{o}s-R\'{e}nyi random graph at criticality. More precisely, the Erd\H{o}s-R\'{e}nyi random graph $G(n,p)$ is obtained by running bond percolation with edge retention probability $p$ on the complete graph with $n$ vertices. In the critical window (that is, when $p=n^{-1}+\lambda n^{-4/3}$ for some $\lambda\in\mathbb{R}$), it was shown in \cite{ABG} that the graphs could be rescaled to a limiting metric space consisting of a tree `glued' together at a finite number of pairs of points. This picture was applied in \cite{CRG} to describe the scaling of the associated random walks, with the limit being described in terms of a `fused' resistance form. In this section, we show that the operation of fusing resistance forms at disjoint pairs of points is continuous with respect to the Gromov-Hausdorff-vague topology (see Proposition \ref{fusingprop}). Consequently, the results of this article give a general strategy for establishing results such as that proved in \cite{CRG}.

Let us start by introducing the operation of fusing an element $(F,R,\mu,\rho)\in\mathbb{F}_c$ over a finite number of subsets. (The restriction to compact spaces is to keep the presentation brief.) Suppose $(V_j)_{j=1}^J$ are a collection of non-empty disjoint compact subsets of $F$, and write
\[\tilde{F}:=\left(F\backslash\bigcup_{j=1}^JV_j\right)\cup\left(\bigcup_{j=1}^J\{V_j\}\right),\]
i.e.\ we consider each of the subsets as a single point. Let $\pi:F\rightarrow\tilde{F}$ be the canonical projection map that is the identity on $F\backslash\cup_{j=1}^JV_j$, and maps points $x\in V_j$ to $\{V_j\}$, $j=1,\dots,J$. Then, define
\begin{equation}\label{tilde1}
\tilde{\mathcal{E}}(f,f):=\mathcal{E}(f\circ\pi,f\circ\pi),\qquad \forall f\in\tilde{\mathcal{F}},
\end{equation}
where
\begin{equation}\label{tilde2}
\tilde{\mathcal{F}}:=\left\{f:\:f\circ\pi\in\mathcal{F}\right\}=\left\{f\circ\pi^{-1}:\:f\in\mathcal{F},\:f|_{V_j}\mbox{ constant, }\forall j\in\{1,\dots,J\}\right\}.
\end{equation}
The form $(\tilde{\mathcal{E}},\tilde{\mathcal{F}})$ is the fused form, and the next lemma ensures it is indeed a resistance form. In the following result, we write $\tilde{\mu}:=\mu\circ\pi^{-1}$, $\tilde{\rho}:=\pi(\rho)$.

{\lem[{cf.\ \cite[Theorem 4.3]{Kig} and \cite[Proposition 2.3]{CRG}}] In the setting described above, $(\tilde{\mathcal{E}},\tilde{\mathcal{F}})$ is a regular resistance form on $\tilde{F}$. Moreover, if $\tilde{R}$ is the associated resistance metric, then $(\tilde{F},\tilde{R},\tilde{\mu},\tilde{\rho})\in\mathbb{F}_c$.}
\begin{proof} Checking RF1 and RF5 of Definition \ref{resformdef} for $(\tilde{\mathcal{E}},\tilde{\mathcal{F}})$ is straightforward. RF2 is also not difficult, though let us explain how to check the completeness of the domain. Suppose $(f_n)_{n\geq 1}$ is Cauchy with respect to $\tilde{\mathcal{E}}$. Then $(f_n\circ \pi)_{n\geq 1}$ is Cauchy with respect to $\mathcal{E}$, so has a limit $g\in\mathcal{F}$. If $u,v\in V_j$ for some $j$, then $f_n\circ\pi(u)=f_n\circ\pi(v)$, and so
\[\left|g(u)-g(v)\right|=\left|g(u)-f_n\circ\pi(u)-g(v)+f_n\circ\pi(v)\right|\leq R(u,v)\mathcal{E}\left(g-f_n\circ\pi,g-f_n\circ\pi\right)\rightarrow 0.\]
Hence $f:=g\circ\pi^{-1}\in\tilde{\mathcal{F}}$, and $f_n\rightarrow f$ with respect to $\tilde{\mathcal{E}}$. We next check RF3. Let $x,y\in \tilde{F}$, $x\neq y$, and set $U=\pi^{-1}(x)$, $V=\pi^{-1}(y)$. Moreover, write $W=(V\cup\cup_{j=1}^JV_j)\backslash U$. Note that $U$ and $W$ are disjoint compact subsets of $F$. Now consider $g_W(u,\cdot)$ for $u\in U$, where we are applying the definition of the resolvent kernel from Lemma \ref{kernellem}. From the latter result, we deduce that $g_W(u,\cdot)|_W=0$, and $g_W(u,u)=R(u,W)>0$ (where the strict positivity is a consequence of \cite[Theorem 4.1 and Corollary 5.4]{Kig}). Since functions in $\mathcal{F}$ are continuous, for each $u\in U$, there exists a $\varepsilon_u>$ such that $g_W(u,\cdot)|_{B_R(u,\varepsilon_u)}>0$. Thus a compactness argument yields the existence of a finite sequence $(u_k)_{k=1}^K$ in $U$ such that $U\subseteq \cup_{k=1,\dots,K}B_R(u_k,\varepsilon_{u_k})$. Write now
\[g:=\sum_{k=1}^Kg_W(u_k,\cdot),\]
so that, by construction, $g|_W=0$ and $g>0$ on $U$. Continuity and compactness implies moreover that there exists an $\varepsilon>0$ such that $g>\varepsilon$ on $U$, and so $\tilde{g}:=\min\{g,\varepsilon\}\in\mathcal{F}$ is constant on $U$. Hence setting $f:=\tilde{g}\circ\pi^{-1}$ gives a function in $\tilde{\mathcal{F}}$ for which $f(x)\neq f(y)$. For RF4, we simply observe that
\begin{equation}\label{rrbound}
\tilde{R}\left(\pi(x),\pi(y)\right)\leq R(x,y),\qquad \forall x,y\in F.
\end{equation}
This completes the proof that $(\tilde{\mathcal{E}},\tilde{\mathcal{F}})$ is a resistance form on $\tilde{F}$. Moreover, (\ref{rrbound}) yields that $(\tilde{F},\tilde{R})$ is compact. Regularity of $(\tilde{\mathcal{E}},\tilde{\mathcal{F}})$ therefore follows from \cite[Corollary 6.4]{Kig}. Checking the remaining requirements for $(\tilde{F},\tilde{R},\tilde{\mu},\tilde{\rho})\in\mathbb{F}_c$ is easy.
\end{proof}

We next specialise to the case when each of the subsets $(V_j)_{j=1}^J$ are given by pairs of points. In particular, we suppose $V_j=\{u_j,v_j\}$ for each $j=1,\dots,J$, so that the elements of interest can be written $(F,R,\mu,\rho,u_1,v_1,\dots,u_J,v_J)$, where $(F,R,\mu,\rho)\in\mathbb{F}_c$. For convenience, we assume henceforth that $J$ is fixed. We define a distance $\Delta_c'$ between $(F,R,\mu,\rho,u_1,v_1,\dots,u_J,v_J)$ and $(F',R',\mu',\rho',u_1',v_1',\dots,u_J',v_J')$ by adding $\sum_{j=1}^J(d_M(\psi(u_j),\psi'(u_j'))+d_M(\psi(v_j),\psi'(v_j')))$ inside the infimum at (\ref{deltacdef}). It is possible to check via standard arguments that this yields a metric (cf.\ \cite[Section 6]{Mgen}), and we call the associated topology the Gromov-Hausdorff-vague topology with additional points. The next result gives the main result of this section, which is the continuity of the fusing map with respect to this topology.

{\propn \label{fusingprop} If
\begin{equation}\label{firstconvergence}
\left({F}_n,R_{n},\mu_{n},\rho_{n},u_1^n,v_1^n,\dots,u_J^n,v_J^n\right)\rightarrow
\left({F},R,\mu,\rho,u_1,v_1,\dots,u_J,v_J\right),
\end{equation}
in the Gromov-Hausdorff-vague topology with additional points, where
$({F}_n,R_{n},\mu_n,\rho_n)$, $n\geq 1$, and $(F,R,\mu,\rho)$ are elements of $\mathbb{F}_c$, and $u_1,v_1,\dots,u_J,v_J\in F$ are distinct, then the associated fused spaces satisfy
\begin{equation}\label{secondconvergence}
\left(\tilde{F}_n,\tilde{R}_n,\tilde{\mu}_n,\tilde{\rho}_n\right)\rightarrow\left(\tilde{F},\tilde{R},\tilde{\mu},\tilde{\rho}\right)
\end{equation}
in the Gromov-Hausdorff-vague topology.}
\begin{proof} We start by observing that, since $({F}_n,R_{n},\mu_{n},\rho_{n})\rightarrow({F},R,\mu,\rho)$, then
\[\liminf_{n\rightarrow\infty}\inf_{x\in F_n}\mu_n\left(\bar{B}_n(x,\delta)\right)>0,\qquad \forall \delta>0,\]
see \cite[Theorem 6.1]{ALWGap}. From \eqref{rrbound}, we thus obtain that
\[\liminf_{n\rightarrow\infty}\inf_{x\in \tilde{F}_n}\tilde{\mu}_n\left(\bar{B}_{\tilde{R}_n}(x,\delta)\right)>0,\qquad \forall \delta>0.\]
Hence, by applying \cite[Theorem 6.1]{ALWGap} again, to complete the proof it will suffice to show (\ref{secondconvergence}) in the (pointed) Gromov-weak topology (that is, for the topology generated by distance given by \eqref{deltacdef}, but without the $d_M^H(\psi(F),\psi'(F))$ term). Applying the alternative definition of the latter topology from \cite[Definition 2.5]{ALW}, this is equivalent to checking
\begin{equation}\label{limittocheck}
\int f\left((\tilde{R}_n(x_i,x_j))_{0\leq i,j\leq K}\right) \prod_{k=1}^K\tilde{\mu}_n(dx_k)
\rightarrow \int f\left((\tilde{R}(x_i,x_j))_{0\leq i,j\leq K}\right) \prod_{k=1}^K\tilde{\mu}(dx_k),
\end{equation}
for any continuous bounded function $f:\mathbb{R}^{K(K+1)/2}\rightarrow\mathbb{R}$, where $x_0=\rho_n$ on the left-hand side, and $x_0=\rho$ on the right-hand side above. Now, note that we can write
\[\int f\left((\tilde{R}(x_i,x_j))_{0\leq i,j\leq K}\right) \prod_{k=1}^K\tilde{\mu}(dx_k)=\int f\left(({R}(\pi(x_i),\pi(x_j)))_{0\leq i,j\leq K}\right) \prod_{k=1}^K{\mu}(dx_k),\]
and similarly for the left-hand side of (\ref{limittocheck}). Moreover, we have that the convergence at (\ref{firstconvergence}) yields convergence of the relevant objects suitably embedded into a common compact metric space $(M,d_M)$, as in the conclusion of Lemma \ref{embeddings}, but also with convergence of the additional points. (We shall assume that the objects are embedded in such a way for the remainder of the proof.) As a consequence of these observations, to check (\ref{limittocheck}) it will be sufficient to show
\begin{equation}\label{lim2}
\tilde{R}_n(\pi_n(x_n),\pi_n(y_n))\rightarrow\tilde{R}(\pi(x),\pi(y))
\end{equation}
whenever $d_M(x_n,x)\rightarrow0$ and $d_M(y_n,y)\rightarrow0$, where we write $\pi_n$ for the canonical projection from $F_n$ to $\tilde{F}_n$.

Now, let $x,y\in F$ be such that $x,y,u_1,v_1,\dots,u_J,v_J$ are distinct, and write $V$ for this collection of vertices. Let $\nu$ be the counting measure on $V$ (placing mass one on each vertex), so that $(V,R|_{V\times V},\nu,x)\in\mathbb{F}_c$ (by Lemma \ref{restriction}). Since the state space is discrete, we can write the associated resistance form as
\begin{equation}\label{onv}
\mathcal{E}_V(f,f)=\frac12\sum_{z,w\in V}c_V(z,w)(f(y)-f(x))^2,
\end{equation}
where the conductances $(c_V(z,w))_{z,w\in V}$ are uniquely determined by the resistance $R|_{V\times V}$ \cite[Theorem 1.7]{Kigdendrite}; the domain of $\mathcal{E}_V$ is all functions from $V$ to $\mathbb{R}$. Similarly, $(\pi(V),R|_{\pi(V)\times \pi(V)},\nu\circ\pi^{-1},\pi(x))\in\mathbb{F}_c$, and we have that the associated form can be written
\begin{equation}\label{onpiv}
\mathcal{E}_{\pi(V)}(f,f)=\frac12\sum_{z,w\in \pi(V)}c_{\pi(V)}(z,w)(f(z)-f(w))^2,
\end{equation}
for conductances $(c_{\pi(V)}(z,w))_{z,w\in \pi(V)}$, and the associated domain is all functions from $\pi(V)$ to $\mathbb{R}$. Now, by \cite[Theorem 8.4]{Kig}, we note that an alternative characterisation of $\mathcal{E}_V$ is as the trace of $(\mathcal{E},\mathcal{F})$ on $V$, namely
\begin{equation}\label{trace1}
\mathcal{E}_V(f,f)=\inf\left\{\mathcal{E}(g,g):\:g\in\mathcal{F},\:g|_{V}=f\right\}.
\end{equation}
Correspondingly, we have that
\begin{equation}\label{trace2}
\mathcal{E}_{\pi(V)}(f,f)=\inf\left\{\tilde{\mathcal{E}}(g,g):\:g\in\tilde{\mathcal{F}},\:g|_{\pi(V)}=f\right\}.
\end{equation}
From (\ref{trace1}), (\ref{trace2}), and the definition of $(\tilde{\mathcal{E}},\tilde{\mathcal{F}})$ at (\ref{tilde1}) and (\ref{tilde2}), it is possible to deduce that $\mathcal{E}_{\pi(V)}(f,f)=\mathcal{E}_V(f\circ\pi,f\circ\pi)$. Putting this together with (\ref{onv}) and (\ref{onpiv}), it follows that
\begin{equation}\label{condsum}
c_{\pi(V)}(z,w)=\sum_{\substack{z',w'\in V:\\\pi(z')=z,\:\pi(w')=w}}c_V(z',w')
\end{equation}
(cf.\ the argument of \cite[Lemma 2.2]{CRG}). We can proceed in exactly the same way for elements of an approximating sequence. Indeed, writing $V_n$ to be the subset of $F_n$ consisting of $x_n,y_n,u_1^n,v_1^n,\dots,u_J^n,v_J^n$, where $d_M(x_n,x)\rightarrow0$ and $d_M(y_n,y)\rightarrow0$, we obtain conductances $(c_{V_n}(z,w))_{z,w\in V_n}$ and $(c_{\pi_n(V_n)}(z,w))_{z,w\in \pi_n(V_n)}$. Moreover, since for large $n$ it is the case that $V_n$ has the same cardinality as $V$, we have a natural correspondence between the latter quantities and the conductances on $V$ and $\pi(V)$ respectively. Clearly, by assumption, we have that $(R_n(z,w))_{z,w\in V_n}\rightarrow (R(z,w))_{z,w\in V}$ (using the natural bijection between elements of $V_n$ and $V$, which exists for large $n$). By the argument of \cite[Lemma 3.3]{CHK}, it follows that $(c_{V_n}(z,w))_{z,w\in V_n}\rightarrow (c_V(z,w))_{z,w\in V}$. Hence, by \eqref{condsum} and the analogous formula for the $n$th level objects, we find that $(c_{\pi_n(V_n)}(z,w))_{z,w\in \pi_n(V_n)}\rightarrow (c_{\pi(V)}(z,w))_{z,w\in \pi(V)}$, from which it readily follows that $(\tilde{R}_n(z,w))_{z,w\in \pi_n(V_n)}\rightarrow (\tilde{R}(z,w))_{z,w\in \pi(V)}$. In particular, this confirms (\ref{lim2}) when $x,y,u_1,v_1,\dots,u_J,v_J$ are distinct.

Next, consider the case when $x,u_1,v_1,\dots,u_J,v_J$ are distinct, and $y=u_1$. Again suppose that $x_n,y_n\in F_n$ are such that $d_M(x_n,x)\rightarrow0$ and $d_M(y_n,y)\rightarrow0$. Then define $V:=\{x,u_1,v_1,\dots,u_J,v_J\}$, $V_n:=\{x_n,u_1^n,v_1^n,\dots,u_J^n,v_J^n\}$, and proceed as above to deduce that $\tilde{R}_n(\pi_n(x_n),\pi_n(u_1^n))\rightarrow \tilde{R}(\pi(x),\pi(u_1))$. Since it also holds that $\tilde{R}_n(\pi(y_n),\pi(u_1^n))\leq R_n(y_n,u_1^n)\leq d_M(y_n,y)+d_M(u_1,u_1^n)\rightarrow 0$, we obtain (\ref{lim2}) in this case too. The remaining cases are similar.
\end{proof}

As a corollary of this result and Theorem \ref{mainthm}, we deduce that if we have convergence of spaces and marked points as at (\ref{firstconvergence}), then we also have convergence of processes associated with the fused spaces. Although we will not present details here, \cite{ABG} (see also the summary of results from that article in \cite{CRG}) gives the required assumption for the critical Erd\H{o}s-R\'{e}nyi random graph, from which the above result yields convergence of the associated rescaled simple random walks. Although the intuition regarding fusing was present in \cite{CRG}, the approach there was based on a sample path argument.

\subsection{Time-changed processes}

In \cite{CHK}, time-changes of stochastic processes on resistance forms were studied, with the main illustrative examples being the Liouville Brownian motion (LBM), the Bouchaud trap model (BTM), and the random conductance model (RCM). Since these were comprehensively covered in the latter article, here we will simply highlight where the present results allow extensions beyond the case of uniform volume doubling (UVD) to be made. Firstly, for the LBM, one advantage of the UVD assumption is that it yields strong estimates on the continuity of the Gaussian field defining the Liouville measure in the model. Whilst we do require some regularity for this \cite{Tal}, we do not need as strong a condition as UVD. Indeed, for establishing convergence of LBM on critical Galton-Watson trees to LBM on the continuum random tree, the desired convergence of Gaussian fields follows from results of \cite{JM} in the case of an offspring distribution with exponential tails, even though this model will not satisfy UVD. Secondly, for BTM, Theorem \ref{mainthm} shows that the requirement of UVD made in \cite[Assumption 5.1]{CHK}, specifically through the appeal to \cite[Assumption 1.2]{CHK}, is not needed; it is possible to obtain convergence of the BTM with the latter replaced by our weaker Assumption \ref{mainassu}. Finally, for the RCM, the argument for trees of \cite[Section 6.1]{CHK} does not depend on UVD, and so \cite[Assumption 1.2]{CHK} can again be replaced by our weaker Assumption \ref{mainassu}. For fractals, however, all the models studied in \cite[Section 6.2]{CHK} naturally satisfy UVD, and so the present work does not extend these results. This is, of course, not to say that the UVD is essential for obtaining scaling limits of stochastic processes on fractals in general; as an example for which we have Gromov-Hausdorff-vague convergence of fractal spaces, but not uniform volume doubling, see \cite{Weiye}.

\bibliography{fractaltrap}

\providecommand{\bysame}{\leavevmode\hbox to3em{\hrulefill}\thinspace}
\providecommand{\MR}{\relax\ifhmode\unskip\space\fi MR }
\providecommand{\MRhref}[2]{%
  \href{http://www.ams.org/mathscinet-getitem?mr=#1}{#2}
}
\providecommand{\href}[2]{#2}
\begin{thebibliography}{10}

\bibitem{ABG}
L.~Addario-Berry, N.~Broutin, and C.~Goldschmidt, \emph{The continuum limit of
  critical random graphs}, Probab. Theory Related Fields \textbf{152} (2012),
  no.~3-4, 367--406.

\bibitem{ABGM}
L.~Addario-Berry, N.~Broutin, C.~Goldschmidt, and G.~Miermont, \emph{The
  scaling limit of the minimum spanning tree of the complete graph}, Ann.
  Probab., to appear.

\bibitem{AEW}
S.~Athreya, M.~Eckhoff, and A.~Winter, \emph{Brownian motion on
  {$\mathbb{R}$}-trees}, Trans. Amer. Math. Soc. \textbf{365} (2013), no.~6,
  3115--3150.

\bibitem{ALW}
S.~Athreya, W.~L\"{o}hr, and A.~Winter, \emph{Invariance principle for variable
  speed random walks on trees}, Ann. Probab., to appear.

\bibitem{ALWGap}
S.~Athreya, W.~L{\"o}hr, and A.~Winter, \emph{The gap between {G}romov-vague
  and {G}romov--{H}ausdorff-vague topology}, Stochastic Process. Appl.
  \textbf{126} (2016), no.~9, 2527--2553.

\bibitem{BCK}
M.~T. Barlow, D.~A. Croydon, and T.~Kumagai, \emph{Subsequential scaling limits
  of simple random walk on the two-dimensional uniform spanning tree}, Ann.
  Probab., to appear.

\bibitem{bACF2}
G.~Ben~Arous, M.~Cabezas, and A.~Fribergh, \emph{Scaling limit for the ant in a
  simple labyrinth}, preprint available at arXiv:1609.03980.

\bibitem{bACF1}
\bysame, \emph{Scaling limit for the ant in high-dimensional labyrinths},
  preprint available at arXiv:1609.03977.

\bibitem{bACF0}
\bysame, \emph{Scaling limit of the random walk on the incipient infinite
  cluster on trees projected to the backbone}, preprint.

\bibitem{Trace}
I.~Benjamini, O.~Gurel-Gurevich, and R.~Lyons, \emph{Recurrence of random walk
  traces}, Ann. Probab. \textbf{35} (2007), no.~2, 732--738.

\bibitem{BBI}
D.~Burago, Y.~Burago, and S.~Ivanov, \emph{A course in metric geometry},
  Graduate Studies in Mathematics, vol.~33, American Mathematical Society,
  Providence, RI, 2001.

\bibitem{Crofin}
D.~A. Croydon, \emph{Convergence of simple random walks on random discrete
  trees to {B}rownian motion on the continuum random tree}, Ann. Inst. Henri
  Poincar\'e Probab. Stat. \textbf{44} (2008), no.~6, 987--1019.

\bibitem{CroyHd}
\bysame, \emph{Hausdorff measure of arcs and {B}rownian motion on {B}rownian
  spatial trees}, Ann. Probab. \textbf{37} (2009), no.~3, 946--978.

\bibitem{RWRRW}
\bysame, \emph{Random walk on the range of random walk}, J. Stat. Phys.
  \textbf{136} (2009), no.~2, 349--372.

\bibitem{Croyrogt}
\bysame, \emph{Scaling limits for simple random walks on random ordered graph
  trees}, Adv. in Appl. Probab. \textbf{42} (2010), no.~2, 528--558.

\bibitem{CRG}
\bysame, \emph{Scaling limit for the random walk on the largest connected
  component of the critical random graph}, Publ. Res. Inst. Math. Sci.
  \textbf{48} (2012), no.~2, 279--338.

\bibitem{CHK}
D.~A. Croydon, B.~M. Hambly, and T.~Kumagai, \emph{Time-changes of stochastic
  processes associated with resistance forms}, preprint available at
  arXiv:1609.02120.

\bibitem{CHKmix}
\bysame, \emph{Convergence of mixing times for sequences of random walks on
  finite graphs}, Electron. J. Probab. \textbf{17} (2012), no. 3, 32.

\bibitem{DLeG}
T.~Duquesne and J.-F. Le~Gall, \emph{Probabilistic and fractal aspects of
  {L}\'evy trees}, Probab. Theory Related Fields \textbf{131} (2005), no.~4,
  553--603.

\bibitem{FOT}
M.~Fukushima, Y.~Oshima, and M.~Takeda, \emph{Dirichlet forms and symmetric
  {M}arkov processes}, extended ed., de Gruyter Studies in Mathematics,
  vol.~19, Walter de Gruyter \& Co., Berlin, 2011.

\bibitem{Georg}
G.~Georganopoulos, \emph{Sur l'approximation des fonctions continues par des
  fonctions lipschitziennes}, C. R. Acad. Sci. Paris S\'er. A-B \textbf{264}
  (1967), A319--A321.

\bibitem{GSu}
A.~L. Gibbs and F.~E. Su, \emph{On choosing and bounding probability metrics},
  Int. Stat. Rev. \textbf{70} (2002), no.~3, 419--435.

\bibitem{GMS}
N.~Gigli, A.~Mondino, and G.~Savar{\'e}, \emph{Convergence of pointed
  non-compact metric measure spaces and stability of {R}icci curvature bounds
  and heat flows}, Proc. Lond. Math. Soc. (3) \textbf{111} (2015), no.~5,
  1071--1129.

\bibitem{HM}
B.~Haas and G.~Miermont, \emph{Scaling limits of {M}arkov branching trees with
  applications to {G}alton-{W}atson and random unordered trees}, Ann. Probab.
  \textbf{40} (2012), no.~6, 2589--2666.

\bibitem{Weiye}
B.~M. Hambly and W.~Yang, \emph{Degenerate limits of non-fixed-point diffusions
  on fractals}, in preparation.

\bibitem{HS}
T.~Hara and G.~Slade, \emph{The scaling limit of the incipient infinite cluster
  in high-dimensional percolation. {II}. {I}ntegrated super-{B}rownian
  excursion}, J. Math. Phys. \textbf{41} (2000), no.~3, 1244--1293,
  Probabilistic techniques in equilibrium and nonequilibrium statistical
  physics.

\bibitem{HHH}
M.~Heydenreich, R.~van~der Hofstad, and T.~Hulshof, \emph{Random walk on the
  high-dimensional {IIC}}, Comm. Math. Phys. \textbf{329} (2014), no.~1,
  57--115.

\bibitem{JM}
S.~Janson and J.-F. Marckert, \emph{Convergence of discrete snakes}, J.
  Theoret. Probab. \textbf{18} (2005), no.~3, 615--647.

\bibitem{KL}
N.~Kajino, \emph{{N}eumann and {D}irichlet heat kernel estimates in inner
  uniform domains for local resistance forms}, in preparation.

\bibitem{Kall}
O.~Kallenberg, \emph{Foundations of modern probability}, second ed.,
  Probability and its Applications (New York), Springer-Verlag, New York, 2002.

\bibitem{Kesten}
H.~Kesten, \emph{Subdiffusive behavior of random walk on a random cluster},
  Ann. Inst. H. Poincar\'e Probab. Statist. \textbf{22} (1986), no.~4,
  425--487.

\bibitem{Kigdendrite}
J.~Kigami, \emph{Harmonic calculus on limits of networks and its application to
  dendrites}, J. Funct. Anal. \textbf{128} (1995), no.~1, 48--86.

\bibitem{kig1}
\bysame, \emph{Analysis on fractals}, Cambridge Tracts in Mathematics, vol.
  143, Cambridge University Press, Cambridge, 2001.

\bibitem{Kig}
\bysame, \emph{Resistance forms, quasisymmetric maps and heat kernel
  estimates}, Mem. Amer. Math. Soc. \textbf{216} (2012), no.~1015, vi+132.

\bibitem{Kum}
T.~Kumagai, \emph{Heat kernel estimates and parabolic {H}arnack inequalities on
  graphs and resistance forms}, Publ. Res. Inst. Math. Sci. \textbf{40} (2004),
  no.~3, 793--818.

\bibitem{MR}
M.~B. Marcus and J.~Rosen, \emph{Markov processes, {G}aussian processes, and
  local times}, Cambridge Studies in Advanced Mathematics, vol. 100, Cambridge
  University Press, Cambridge, 2006.

\bibitem{Mgen}
G.~Miermont, \emph{Tessellations of random maps of arbitrary genus}, Ann. Sci.
  \'Ec. Norm. Sup\'er. (4) \textbf{42} (2009), no.~5, 725--781.

\bibitem{Stone}
C.~Stone, \emph{Limit theorems for random walks, birth and death processes, and
  diffusion processes}, Illinois J. Math. \textbf{7} (1963), 638--660.

\bibitem{Suz2}
K.~Suzuki, \emph{Convergence of {B}rownian motions on {RCD}*({K},$\infty$)
  spaces}, preprint available at arXiv:1603.08622.

\bibitem{Suz1}
\bysame, \emph{Convergence of {B}rownian motions on {RCD}*({K},{N}) spaces},
  preprint available at arXiv:1509.02025.

\bibitem{Tal}
M.~Talagrand, \emph{Regularity of {G}aussian processes}, Acta Math.
  \textbf{159} (1987), no.~1-2, 99--149.

\bibitem{HJ}
R.~van~der Hofstad and A.~A. J{\'a}rai, \emph{The incipient infinite cluster
  for high-dimensional unoriented percolation}, J. Statist. Phys. \textbf{114}
  (2004), no.~3-4, 625--663.

\end{thebibliography}
\bibliographystyle{amsplain}

\end{document}